\documentclass[a4paper]{amsart}
\usepackage[british]{babel}

\usepackage[all]{xy}
\usepackage{amssymb,enumitem,nicefrac,mathtools}

\usepackage[pdftitle={Equivariant representable K-theory},
  pdfauthor={Heath Emerson and Ralf Meyer},
  pdfsubject={Mathematics; MSC 19K35, 46L80}
]{hyperref}
\newcommand*{\MRref}[2]{\href{http://www.ams.org/mathscinet-getitem?mr=#1}{MR \textbf{#1}}}
\newcommand*{\arxiv}[1]{\href{http://www.arxiv.org/abs/#1}{arXiv: #1}}

\usepackage[lite]{amsrefs}

\usepackage{microtype}

\numberwithin{equation}{section}

\newtheorem{theorem}[equation]{Theorem}
\newtheorem{lemma}[equation]{Lemma}
\newtheorem{proposition}[equation]{Proposition}
\newtheorem{corollary}[equation]{Corollary}
\theoremstyle{definition}
\newtheorem{definition}[equation]{Definition}
\newtheorem{notation}[equation]{Notation}
\theoremstyle{remark}
\newtheorem{remark}[equation]{Remark}
\newtheorem{example}[equation]{Example}

\newcommand*{\dd}{\textup{d}} % differential forms
\newcommand*{\nc}{_\textup{n}} % subscript for norm-continuous
\newcommand*{\ns}{_\textup{s}} % subscript for *-strongly continuous
\newcommand*{\bound}{_\textup{b}} % bounded elements

\newcommand*{\ID}{\textup{id}}

\DeclareMathOperator{\KK}{KK}
\DeclareMathOperator{\cRKK}{\mathcal{R}KK}
\DeclareMathOperator{\RKK}{RKK}
\DeclareMathOperator{\K}{K}
\DeclareMathOperator{\RK}{RK}
\DeclareMathOperator{\vbK}{VK}
\DeclareMathOperator{\Vect}{Vect}

\DeclareMathOperator{\tsr}{\tsr}

\newcommand*{\Cat}{\mathfrak{C}}
\newcommand*{\Fred}{\mathcal{F}}
\newcommand*{\Unitary}{\mathcal{U}}
\newcommand*{\Prid}{\mathcal{P}}
\newcommand*{\Mult}{\mathcal{M}}

\newcommand*{\Hilm}{\mathcal{E}}
\newcommand*{\HilmF}{\mathcal{F}}
\newcommand*{\Hils}{\mathcal{H}}
\newcommand*{\CONT}{\mathcal{C}}

\newcommand*{\Grd}{\mathcal{G}}
\newcommand*{\Base}{Z}
\newcommand*{\base}{z}
\newcommand*{\EG}{\mathcal{E}}

\newcommand*{\C}{\mathbb C}
\newcommand*{\Z}{\mathbb Z}
\newcommand*{\R}{\mathbb R}
\newcommand*{\N}{\mathbb N}
\newcommand*{\Sphere}{\mathbb S}

\newcommand*{\Mat}{\mathbb M}
\newcommand*{\Comp}{\mathbb K}
\newcommand*{\Bound}{\mathbb B}

\newcommand*{\nb}{\nobreakdash}
\newcommand*{\Cstar}{\texorpdfstring{$C^*$\nobreakdash-}{C*-}}
\newcommand*{\sigCstar}{\texorpdfstring{$\sigma$\nobreakdash-$C^*$-\hspace{0pt}}{sigma-C*-}}
\newcommand*{\Star}{\texorpdfstring{$^*$\nobreakdash-\hspace{0pt}}{*-}}
\newcommand*{\sC}{\sigma\text{-}C^*}

\newcommand*{\abs}[1]{\lvert#1\rvert}
\newcommand*{\norm}[1]{\lVert#1\rVert}
\newcommand*{\cl}[1]{\overline{#1}}
\newcommand*{\conj}[1]{\overline{#1}}

\newcommand*{\defeq}{\mathrel{\vcentcolon=}}
\newcommand*{\into}{\rightarrowtail}
\newcommand*{\prto}{\twoheadrightarrow}

\begin{document}

\title{Equivariant representable K-theory}
%\author{Heath Emerson and Ralf Meyer}
\author{Heath Emerson}
\email{hemerson@math.uvic.ca}

\address{Department of Mathematics and Statistics\\
  University of Victoria\\
  PO BOX 3045 STN CSC\\
  Victoria, B.C.\\
  Canada V8W 3P4}

\author{Ralf Meyer}
\email{rameyer@uni-math.gwdg.de}

\address{Mathematisches Institut\\
  Georg-August Universit\"at G\"ottingen\\
  Bunsenstra{\ss}e 3--5\\
  37073 G\"ottingen\\
  Germany}


\begin{abstract}
  We interpret certain equivariant Kasparov groups as
  equivariant representable \(\K\)\nb-theory groups and compute
  these via a classifying space and as \(\K\)\nb-theory groups
  of suitable \sigCstar{}algebras.  We also relate equivariant
  vector bundles to these \sigCstar{}algebras and provide
  sufficient conditions for equivariant vector bundles to
  generate representable \(\K\)\nb-theory.  We mostly work in
  the generality of locally compact groupoids with Haar system.
\end{abstract}
\maketitle

\section{Introduction}
\label{sec:intro}

For a locally compact space~\(X\), we must distinguish its
\(\K\)\nb-theory \(\K^*(X)\) and its representable
\(\K\)\nb-theory \(\RK^*(X)\).  Vector bundles on~\(X\) define
classes in \(\RK^*(X)\), but not necessarily in~\(\K^*(X)\).
Furthermore, given a map \(X\to Y\), we can define
\(\RK^*_Y(X)\), the \(\K\)\nb-theory of~\(X\) with
\(Y\)\nb-compact support.  Such groups are important, for
instance, because the symbol of an elliptic pseudodifferential
operator on a smooth manifold~\(X\) lies in the group
\(\RK^*_X(TX)\).

We may define all three theories using homotopy classes of
continuous maps from~\(X\) to the space of Fredholm operators
on a Hilbert space.  For \(\RK^0(X)\), we allow all such maps,
for \(\K^0(X)\), we require the map to have unitary values
outside a compact subset, and for \(\RK^*_Y(X)\), we require
unitary values outside a \(Y\)\nb-compact subset of~\(X\).

Alternatively, we may identify \(\RK^*(X)\) and \(\RK^*_Y(X)\)
with the equivariant Kasparov groups
\(\KK^X_*\bigl(\CONT_0(X),\CONT_0(X)\bigr)\) and
\(\KK^Y_*\bigl(\CONT_0(Y),\CONT_0(X)\bigr)\), respectively,
where we view \(X\) and~\(Y\) as groupoids with only identity
morphisms (these groups are denoted by \(\RKK(X;\C,\C) \cong
\cRKK\bigl(X;\CONT_0(X),\CONT_0(X)\bigr)\) and
\(\cRKK\bigl(Y;\CONT_0(Y),\CONT_0(X)\bigr)\)
in~\cite{Kasparov:Novikov}).

In this article, we study equivariant versions of these three
kinds of \(\K\)\nb-theory.  We establish an equivalence between
several alternative definitions that are based on equivariant
Kasparov theory, maps to spaces of Fredholm operators,
\(\K\)\nb-theory of certain crossed product
\sigCstar{}algebras, and equivariant vector bundles,
respectively.  The main goal here is to show that certain
equivariant Kasparov groups deserve to be called
``representable equivariant \(\K\)\nb-theory'' or ``equivariant
\(\K\)\nb-theory with \(Y\)\nb-compact support.''  Moreover, we
study the general properties of these theories and investigate
when representable equivariant \(\K\)\nb-theory can be
described by equivariant vector bundles.

Although we are mainly interested in the case of group actions,
we must consider crossed product groupoids \(G\ltimes X\),
anyway.  Therefore, we work with groupoids most of the time.
More precisely, we use Hausdorff, second countable, locally
compact groupoids with Haar systems.  Let~\(\Grd\) be such a
groupoid.

Let~\(X\) be a (Hausdorff, locally compact, second countable)
proper \(\Grd\)\nb-space, that is, the crossed product groupoid
\(\Grd\ltimes X\) is proper.  We define the
\(\Grd\)\nb-equivariant \(\K\)\nb-theory of~\(X\) by
\[
\K^*_\Grd(X) \defeq \K_*\bigl(\Grd\ltimes \CONT_0(X)\bigr).
\]
If~\(\Grd\) is a group, this agrees with Chris Phillips'
definition in~\cite{Phillips:Equivariant_Kbook}.

The \(\Grd\)\nb-equivariant representable \(\K\)\nb-theory
of~\(X\) is defined by
\[
\RK^*_\Grd(X) \defeq
\KK_*^{\Grd\ltimes X}\bigl(\CONT_0(X),\CONT_0(X)\bigr).
\]

Let~\(Y\) be another \(\Grd\ltimes X\)-space.  The
\(\Grd\)\nb-equivariant \(\K\)\nb-theory of~\(Y\) with
\(X\)\nb-compact supports is defined by
\[
\RK^*_{\Grd,X}(Y) \defeq
\KK_*^{\Grd\ltimes X}\bigl(\CONT_0(X),\CONT_0(Y)\bigr).
\]

Of course, \(\RK^*_\Grd(X)\) is the special case
\(\RK^*_{\Grd,X}(X)\) of the last definition.  Conversely,
\(\RK^*_{\Grd,X}(Y)\) can be computed as the inductive limit of
the groups \(\RK^*_\Grd(A,\partial A)\), where~\(A\) runs
through the directed set of \(X\)\nb-compact
\(\Grd\)\nb-invariant subsets of~\(Y\) and
\(\RK^*_\Grd(A,\partial A)\) denotes the relative version of
\(\RK^*_\Grd\) (see Theorem~\ref{the:RK_support_continuous}).

A first alternative definition of \(\RK^*_\Grd(X)\) extends
Graeme Segal's description of representable \(\K\)\nb-theory as
the set of homotopy classes of maps to the space of Fredholm
operators on a separable Hilbert space
(\cite{Segal:Fredholm_complexes}).  We use the direct sum of
countably many copies of the regular representation,
\(\Hils_\Grd \defeq \bigoplus_{n\in\N} L^2(\Grd)\).  This is a
continuous field of Hilbert spaces over the object
space~\(\Grd^{(0)}\) of the groupoid.  Let~\(\Fred_\Grd\) be
the set over~\(\Grd^{(0)}\) of essentially unitary operators on
the corresponding fibres of~\(\Hils_\Grd\).  The
groupoid~\(\Grd\) acts on~\(\Fred_\Grd\) in an obvious way.
Almost by definition, the Kasparov group \(\RK^0_\Grd(X)\) is
the space of homotopy classes of \(\Grd\)\nb-equivariant maps
\(X\to\Fred_\Grd\) that are continuous in a suitable sense.
But we must be careful about the topology here.

In the groupoid case, we cannot expect the continuous field of
Hilbert spaces~\(\Hils_\Grd\) to be locally trivial.  As a
result, there is no canonical norm topology on~\(\Fred_\Grd\);
only on compact operators is the norm topology canonically
defined.  We equip~\(\Fred_\Grd\) with the topology where a net
\((U_i)\) converges to~\(U\) if and only if \(U_i\to U\) and
\(U_i^*\to U^*\) strongly and \(1-U_iU_i^*\to 1-UU^*\) and
\(1-U_i^*U_i \to 1-U^*U\) in norm.  For this topology,
\(\RK^0_\Grd(X)\) is the space of homotopy classes of
\(\Grd\)\nb-equivariant continuous maps \(X\to\Fred_\Grd\).

There are similar descriptions for \(\RK^0_{\Grd,Y}(X)\) and
\(\K^0_\Grd(Y)\).  The support of a continuous
\(\Grd\)\nb-equivariant map \(U\colon X\to\Fred_\Grd\) is the
closure of the set of \({x\in X}\) where~\(U_x\) is not
unitary.  We identify \(\RK^0_{\Grd,Y}(X)\) and
\(\K^0_\Grd(Y)\) with the sets of homotopy classes of
continuous \(\Grd\)\nb-equivariant maps \(X\to\Fred_\Grd\) with
\(Y\)\nb-compact or \(\Grd\)\nb-compact support, respectively.

For a locally compact group~\(G\), there is a canonical norm
topology on~\(\Fred_G\).  We show that it makes no difference
whether we use this norm topology instead of the weaker
topology described above.  This requires, among other things,
an equivariant version of the Cuntz--Higson--Kuiper Theorem on
the contractibility of the unitary group of stable multiplier
algebras of \(\sigma\)\nb-unital \Cstar{}algebras
(see~\cite{Cuntz-Higson:Kuiper}).

We also describe \(\RK^*_{\Grd,X}(Y)\) as the
\(\K\)\nb-theory of a certain \sigCstar{}algebra
(\cites{Phillips:Inverse, Phillips:Representable_K}), that is,
a countable inverse limit of \Cstar{}algebras.  Let
\((X_n)_{n\in\N}\) be an increasing sequence of
\(\Grd\)\nb-compact \(\Grd\)\nb-invariant subsets of~\(X\)
with \(\bigcup X_n=X\), and let~\(Y_n\) be the pre-image
of~\(X_n\) with respect to the canonical map \(Y\to X\).  Then
\[
\RK^*_{\Grd,X}(Y) \cong
\K_*\bigl(\varprojlim C^*(\Grd\ltimes Y_n)\bigr).
\]
Here we use Phillips' representable \(\K\)\nb-theory for
\sigCstar{}algebras, so that our~\(\K_0\) does not agree with
algebraic \(\K_0\).  Using the exactness properties of
\sigCstar{}algebra \(\K\)\nb-theory, we define
\(\RK^*_\Grd(X,A)\) for pairs and establish excision and long
exact sequences for this theory.  Furthermore, we get a Milnor
\(\varprojlim\nolimits^1\)-sequence that relates
\(\RK_{\Grd,X}^*(Y)\) to the projective system
\(\bigl(\K_\Grd^*(Y_n)\bigr)_{n\in\N}\).  Other properties like
homotopy invariance and induction isomorphisms are immediate
from the definitions as well.

As in~\cite{Lueck-Oliver:Completion}, we also consider a
version of \(\K\)\nb-theory defined by equivariant vector
bundles.  Let \(\vbK^0_\Grd(X)\) be the Grothendieck group of
the monoid of \(\Grd\)\nb-equivariant vector bundles on~\(X\).
We extend this to a \(\Z/2\)-graded theory for pairs of spaces.
There is always a canonical map
\begin{equation}
  \label{eq:intro_compact_vbK_RK}
  \vbK^*_\Grd(X,A)\to \RK^*_\Grd(X,A),
\end{equation}
but it need not be an isomorphism in general, as various
counterexamples in \cites{Phillips:Equivariant_K2,
  Lueck-Oliver:Completion, Sauer:K-theory} show.  We provide a
sufficient condition for an isomorphism, but it is not
particularly striking.  Similar results about twisted
equivariant \(\K\)\nb-theory for Lie groupoids are established
in~\cite{Tu-Xu-Laurent-Gengoux:Twisted_K}.

Our approach is to relate equivariant vector bundles to the same
\sigCstar{}algebra
\[
\sC(\Grd\ltimes X) \defeq \varprojlim C^*(\Grd\ltimes X_n)
\]
that already computes \(\RK^*_\Grd(X)\).  We show that
\(\Grd\)\nb-equivariant vector bundles on~\(X\) correspond to
projections in the stabilisation
\[
\sC(\Grd\ltimes X)_\Comp \defeq
\varprojlim C^*(\Grd\ltimes X_n) \otimes \Comp(\ell^2\N).
\]
Hence the question is whether \(\K_0\bigl(\sC(\Grd\ltimes
X)\bigr)\) is the Grothendieck group of the monoid of
projections in \(\sC(\Grd\ltimes X)_\Comp\).  This happens if
\(\sC(\Grd\ltimes X)_\Comp\) is a \Cstar{}algebra with an
approximate unit of projections.  Notice that \(\sC(\Grd\ltimes
X)_\Comp\) is a \Cstar{}algebra if and only if~\(G\) acts
cocompactly on~\(X\).  Our positive results are limited to
cocompact actions for this reason.

We show that \(\sC(\Grd\ltimes X)_\Comp\) has an approximate
unit of projections if and only if for each \(x\in X\) and each
irreducible representation of the stabiliser of~\(x\), there is
a \(\Grd\)\nb-equivariant vector bundle on~\(X\) whose
restriction to~\(x\) contains the given representation of the
stabiliser.  If, in addition, the action of~\(\Grd\) on~\(X\)
is cocompact as well, then the map
in~\eqref{eq:intro_compact_vbK_RK} is an isomorphism for all
closed \(\Grd\)\nb-invariant subsets \(A\subseteq X\) and all
\(*\in\Z\).

This sufficient criterion covers all cases where we
know~\eqref{eq:intro_compact_vbK_RK} to be an isomorphism:
cocompact actions of almost connected groups, discrete groups,
or inverse limits of discrete groups (see
\cites{Phillips:Equivariant_K2, Lueck-Oliver:Completion,
  Sauer:K-theory}), and some simple cases like equivalence
relations, orbifold groupoids, and group actions that are
induced from compact subgroups.  A counterexample by Juliane
Sauer yields examples of a totally disconnected group or of a
bundle of totally disconnected compact groups over the circle
for which \(\vbK^0_\Grd(X) \not\cong \RK^0_\Grd(X)\) for some
space~\(X\).

Since any proper group action is locally induced from a compact
subgroup, the problem for cocompact group actions is a purely
global one.  The counterexamples in
\cites{Phillips:Equivariant_Kbook, Lueck-Oliver:Completion,
  Sauer:K-theory} and the example of continuous trace
\Cstar{}algebras show that the existence of an approximate unit
of projections in the stabilisation is a rather subtle
property.

We finish the introduction with an application that motivated
us to write this article.  Let \(X\) and~\(Y\) be two smooth
manifolds.  Then a family of elliptic pseudodifferential
operators on~\(X\) parametrised by~\(Y\), on the one hand,
provides an element in \(\KK_*\bigl(\CONT_0(X),
\CONT_0(Y)\bigr)\) and, on the other hand, has a symbol in the
\(\K\)\nb-theory of \(TX\times Y\) with \(X\)\nb-compact
support, where \(TX\) denotes the tangent space of~\(X\).  This
remains true equivariantly.  Let~\(G\) be a locally compact
group, let~\(Y\) be a \(G\)\nb-space, and let~\(X\) be a
proper, smooth \(G\)\nb-manifold with tangent bundle~\(TX\).
The second Poincar\'e duality isomorphism
in~\cite{Emerson-Meyer:Dualities} provides a natural
isomorphism
\[
\KK_*^G\bigl(\CONT_0(X),\CONT_0(Y)\bigr) \cong
\RK^*_{G,X}(TX\times Y).
\]
The proof of this isomorphism in~\cite{Emerson-Meyer:Dualities}
depends on our definition of \(\RK^*_{G,X}(TX\times Y)\) as an
equivariant Kasparov group.  It provides more geometric
descriptions of \(\KK_*^G\bigl(\CONT_0(X),\CONT_0(Y)\bigr)\)
when we combine it with our alternative descriptions of
\(\RK^*_{G,X}(TX\times Y)\) in terms of equivariant maps to a
space of Fredholm operators or, in nice cases, equivariant
vector bundles.

\section{Basic definitions}
\label{sec:definitions}

Throughout this article, we act as if we were dealing with
\emph{complex} \Cstar{}algebras and complex \(\K\)\nb-theory.
But everything carries over to the real case with obvious
modifications.  All spaces shall be locally compact, Hausdorff,
and second countable, and all \Cstar{}algebras separable.  We
work equivariantly with respect to groupoids most of the time.
We usually write~\(\Grd\) for a groupoid and~\(G\) for a group.
All groupoids are --~usually tacitly~-- required locally
compact, Hausdorff, second countable, and with a Haar system
(see \cite{Renault:Groupoid_Cstar}*{Definition 2.2}).  We do
not consider non-Hausdorff groupoids because they have few
proper actions on Hausdorff spaces.

Let~\(\Grd\) be a groupoid (as above).  We write \(\Base\defeq
\Grd^{(0)}\) for its object space, \(\Grd^{(1)}\) for its
morphism space, and \(r,s\colon \Grd^{(1)}\rightrightarrows
\Grd^{(0)}\) for the range and source maps.  By our convention,
the spaces \(\Grd^{(1)}\) and \(\Grd^{(0)}\) are locally
compact, Hausdorff, second countable, and the range, source,
unit, multiplication and inversion maps are continuous.  The
existence of a Haar system forces the range, source, and
multiplication maps to be open (see
\cite{Westman:Non-transitive}
or~\cite{Paterson:Groupoids}*{Proposition 2.2.1}).

Actions of groupoids on spaces and \Cstar{}algebras are
explained, for instance, in~\cite{LeGall:KK_groupoid}, which is
the basic reference for Kasparov theory for groupoids.  We will
use this theory and its main properties throughout, without
much comment.

Recall that a space over~\(\Grd^{(0)}\) is a space with a
continuous map to~\(\Grd^{(0)}\).  A \(\Grd\)\nb-space is a
space over~\(\Grd^{(0)}\) with maps \(\alpha_g\colon X_{s(g)}
\to X_{r(g)}\) for all \(g\in G\), where \(X_{s(g)}\)
and~\(X_{r(g)}\) are the fibres of~\(X\) over \(s(g)\) and
\(r(g)\); the maps~\(\alpha_g\) vary continuously in the sense
that they combine to a continuous map from \(s^*(X) \defeq
\Grd^{(1)} \times_s X\) to \(r^*(X) \defeq \Grd^{(1)} \times_s
X\).

Crossed products for groupoid actions on \Cstar{}algebras are
defined in great generality in~\cite{Renault:Representations}.
We denote the \Cstar{}algebra crossed product by \(\Grd\ltimes
A\) and groupoid \Cstar{}algebras by \(C^*(\Grd)\).  Readers
unfamiliar with this may restrict attention to the well-known
case of group actions and their crossed products.  Since we
only consider proper actions here, full and reduced crossed
products agree in all cases we need.

Given an action of a groupoid~\(\Grd\) on a space~\(X\), we
often form the crossed product groupoid \(\Grd\ltimes X\),
which is another groupoid with object space~\(X\) and morphisms
encoding the action.  We have a canonical \Cstar{}algebra
isomorphism
\begin{equation}
  \label{eq:crossed_product_groupoid}
  C^*(\Grd\ltimes X) \cong \Grd \ltimes \CONT_0(X)
\end{equation}
because both algebras satisfy the same universal property.

Let~\(\Grd\) be a (locally compact, Hausdorff, etc.) groupoid
and~\(X\) a proper \(\Grd\)\nb-space.  Although the following
definitions make sense for non-proper actions as well, we only
consider the special case of proper actions in this article.

\begin{definition}
  \label{def:equivariant_K}
  The \emph{\(\Grd\)\nb-equivariant \(\K\)\nb-theory of~\(X\)}
  is defined by
  \[
  \K_\Grd^*(X) \defeq
  \K_*\bigl(\Grd\ltimes \CONT_0(X)\bigr)
  \cong \K_*\bigl(C^*(\Grd\ltimes X)\bigr);
  \]
  here we use the canonical
  isomorphism~\eqref{eq:crossed_product_groupoid}.
\end{definition}

If~\(\Grd\) is a group, then \(\K_\Grd^*(X)\) agrees with the
equivariant \(\K\)\nb-theory for locally compact groups studied
by Chris Phillips in~\cite{Phillips:Equivariant_Kbook}.  We
will soon spell out a Kasparov type description of
\(\K_\Grd^*(X)\) using Fredholm operators on continuous fields
of Hilbert spaces, which is more similar to Phillips'
definition.  Moreover, Definition~\ref{def:equivariant_K} is a
special case of the definition of twisted equivariant
\(\K\)\nb-theory in~\cite{Tu-Xu-Laurent-Gengoux:Twisted_K}.

\begin{definition}
  \label{def:equivariant_RK}
  The \emph{\(\Grd\)\nb-equivariant representable
    \(\K\)\nb-theory of~\(X\)} is
  \[
  \RK^*_\Grd(X) \defeq
  \KK^{\Grd\ltimes X}_*\bigl(\CONT_0(X),\CONT_0(X)\bigr).
  \]
\end{definition}

\begin{definition}
  \label{def:X-compact}
  Let~\(Y\) be a \(\Grd\ltimes X\)\nb-space via some map
  \(p\colon Y\to X\).  A subset~\(A\) of~\(Y\) is
  \emph{\(X\)\nb-compact} if \(p\colon Y\to X\) restricts to a
  proper map \(A\to X\).
\end{definition}

\begin{definition}
  \label{def:K-support}
  Let~\(Y\) be a \(\Grd\ltimes X\)\nb-space via some map
  \(p\colon Y\to X\).  The \emph{\(\Grd\)\nb-equivariant
    \(\K\)\nb-theory of~\(Y\) with \(X\)\nb-compact support} is
  \[
  \RK_{\Grd,X}^*(Y) \defeq
  \KK^{\Grd\ltimes X}_*\bigl(\CONT_0(X),\CONT_0(Y)\bigr).
  \]
\end{definition}

This definition contains \(\RK_\Grd^*(X) = \RK_{\Grd,X}^*(X)\)
as a special case.  In contrast, \(\K_\Grd^*(X)\) need not be a
special case of Definition~\ref{def:K-support}.  Our
terminology is justified by Theorem~\ref{the:RK_Fred_general},
which describes \(\RK^0_\Grd(X)\), \(\RK^0_{\Grd,X}(Y)\), and
\(\K^0_\Grd(X)\) using spaces of maps that only differ by
support conditions.

\begin{definition}
  \label{def:vbK}
  A \emph{\(\Grd\)\nb-equivariant vector bundle over~\(X\)} is
  a vector bundle~\(V\) over~\(X\) together with a continuous,
  fibrewise linear \(\Grd\)\nb-action on its total space.

  Isomorphism classes of \(\Grd\)\nb-equivariant vector bundles
  over~\(X\) form a monoid \(\Vect_\Grd(X)\) with respect to
  the usual direct sum of vector bundles.  Let
  \(\vbK^0_\Grd(X)\) be the Grothendieck group of this monoid.
\end{definition}

We can also extend Definition~\ref{def:vbK} to pairs of spaces.
By definition, a \emph{pair of \(\Grd\)\nb-spaces} is a
\(\Grd\)\nb-space~\(X\) together with a closed
\(\Grd\)\nb-invariant subspace~\(A\).

\begin{definition}
  \label{def:vbK_relative}
  Let \((X,A)\) be a pair of \(\Grd\)\nb-spaces.  Consider
  triples \((V^+,V^-,\varphi)\) where \(V^\pm\in\Vect_\Grd(X)\)
  and \(\varphi\colon V^+|_A \to V^-|_A\) is a
  \(\Grd\)\nb-equivariant vector bundle isomorphism.
  Isomorphism classes of such triples form a monoid with
  respect to direct sum, which we denote by
  \(\Vect_\Grd(X,A)\).

  Call a triple \((V^+,V^-,\varphi)\) \emph{degenerate}
  if~\(\varphi\) extends to a \(\Grd\)\nb-equivariant vector
  bundle isomorphism \(V^+\to V^-\).  The degenerate triples
  form a submonoid of \(\Vect_\Grd(X,A)\).  Call two triples
  \emph{stably isomorphic} if they become isomorphic after
  adding degenerate triples.  We let \(\vbK^0_\Grd(X,A)\) be
  the monoid of stable isomorphism classes of
  \(\Vect_\Grd(X,A)\).
\end{definition}

By design, the class of degenerate triples is a neutral element
in \(\vbK^0_\Grd(X,A)\).  For any triple \((V^+,V^-,\varphi)\),
the direct sum \((V^+,V^-,\varphi) \oplus
(V^-,V^+,\varphi^{-1})\) is degenerate (compare
\cite{Cuntz-Meyer-Rosenberg}*{Lemma 1.42}).  Hence
\(\vbK^0_\Grd(X,A)\) is a group.  There is a canonical
isomorphism \(\vbK^0_\Grd(X,\emptyset) \cong \vbK^0_\Grd(X)\).

We extend \(\vbK^0_\Grd(X)\) to a graded theory by
\[
\vbK^{-n}_\Grd(X) \defeq
\vbK^0_\Grd(X\times\Sphere^n,X\times\{\star\}),
\]
where~\(\Sphere^n\) denotes the \(n\)\nb-sphere and~\(\star\)
is its base point.  For the relative theory, we put
\begin{equation}
  \label{eq:vbK_odd}
  \vbK^{-n}_\Grd(X,A) \defeq
  \vbK^0_\Grd(X\times\Sphere^n,
  X\times\{\star\}\cup A\times \Sphere^n).
\end{equation}

\begin{remark}
  \label{rem:RK_compact_group}
  Let~\(G\) be a compact group.  The \emph{Green--Julg Theorem}
  (\cite{Julg:K_equivariante}) identifies \(\K^*_G(X)\) with
  the equivariant \(\K\)\nb-theory of~\(X\) studied by
  topologists, and \cite{Kasparov:Novikov}*{Proposition 2.20}
  identifies \(\RK^*_G(X)\) with the \(G\)\nb-equivariant
  representable \(\K\)\nb-theory in the sense of Graeme Segal
  (\cite{Segal:Fredholm_complexes}).  This justifies our
  notation.

  If both \(G\) and~\(X\) are compact, then \(\K^*_G(X) \cong
  \RK^*_G(X) \cong \vbK^*_G(X)\).  If~\(X\) is not compact,
  then \(\K^*_G(X) \cong \vbK^*_G(X^+,\{\infty\})\), where
  \(X^+=X\cup\{\infty\}\) is the one-point compactification
  of~\(X\).  There is a canonical map
  \(\vbK^*_G(X)\to\RK^*_G(X)\) for any~\(X\), which may or may
  not be an isomorphism.
\end{remark}

\section{Description by maps to Fredholm operators}
\label{sec:K_via_Fred}

Our next goal is to describe \(\RK_\Grd^0(X)\) as the space of
homotopy classes of equivariant maps to a suitable space of
Fredholm operators.  Such a description is well-known in the
non-equivariant case, see
\cites{Atiyah:K-theory,Segal:Fredholm_complexes}.  For
equivariant \(\K\)\nb-theory, this is a special case of a
similar result for twisted equivariant \(\K\)\nb-theory
in~\cite{Tu-Xu-Laurent-Gengoux:Twisted_K}.  For general
groupoids, we equip the space of Fredholm operators with an
unusual topology because the norm topology does not make sense.
In~\S\ref{sec:topology} we restrict attention to groups, where
the norm topology makes sense, and show that the result remains
valid if we use the norm topology instead.

\subsection{Regular representation and Equivariant
  Stabilisation Theorem}
\label{sec:regular}

Let~\(\Grd\) be a (locally compact, Hausdorff, etc.) groupoid
and let \(\Base\defeq\Grd^{(0)}\).  For \(\base\in\Base\), let
\[
\Grd_\base\defeq \{g\in\Grd^{(1)} \mid s(g)=\base\}.
\]
The \emph{Haar system} of~\(\Grd\) consists of measures
\((\lambda_\base)_{\base\in\Base}\) such that~\(\lambda_\base\)
has support~\(\Grd_\base\) and
\[
\Base\ni \base\mapsto
\int_{\Grd_\base} f(g)\,\dd\lambda_\base(g)
\]
is a continuous function for all \(f\in\CONT_c(\Grd^{(1)})\)
and \(\dd\lambda_{gz} = g\cdot\dd\lambda_z\) (see
\cite{Renault:Groupoid_Cstar}*{Definition 2.2}).

Define a \(\CONT_0(\Base)\)-valued inner product on
\(\CONT_c(\Grd^{(1)})\) by
\[
\langle f_1,f_2\rangle(\base) \defeq
\int_{\Grd_\base} \conj{f_1(g)}\cdot f_2(g)
\,\dd\lambda_\base(g)
\]
and a right \(\CONT_0(\Base)\)-module structure by pointwise
multiplication.  This yields a pre-Hilbert module, whose
completion is a Hilbert module over \(\CONT_0(\Base)\), which
we denote by \(L^2(\Grd)\).  The usual left regular
representation of~\(\Grd\) on \(\CONT_c(\Grd^{(1)})\) is
compatible with the pre-Hilbert module structure, so that
\(L^2(\Grd)\) becomes a \(\Grd\)\nb-equivariant Hilbert
\(\CONT_0(\Base)\)-module.  This is the \emph{left regular
  representation} of~\(\Grd\).

\begin{definition}
  \label{def:universal_Hilbert_module}
  Let~\(A\) be a \(\Grd\)\nb-\Cstar{}algebra.  Then we define
  \[
  \Hils_{A,\Grd} \defeq
  L^2(\Grd)^\infty \otimes_{\CONT_0(\Base)} A,
  \]
  equipped with the diagonal representation of~\(\Grd\).  Here
  we write \(\Hilm^\infty\) for the Hilbert module direct sum
  of countably many copies of a Hilbert module~\(\Hilm\).
\end{definition}

\begin{theorem}[Equivariant Stabilisation Theorem, see
  \cite{Tu:Novikov}*{Corollaire 6.22}]
  \label{the:stabilisation}
  Let~\(X\) be a proper \(\Grd\)\nb-space, let~\(A\) be a
  \(\Grd\ltimes X\)-\Cstar{}algebra, and let~\(\Hilm\) be a
  countably generated \(\Grd\)\nb-equivariant Hilbert
  \(A\)\nb-module.  Then there is a \(\Grd\)\nb-equivariant
  isomorphism \(\Hilm\oplus \Hils_{A,\Grd} \cong
  \Hils_{A,\Grd}\).
\end{theorem}

This usually allows us to restrict attention to the special
Hilbert module~\(\Hils_{A,\Grd}\).

Hilbert modules over \(\CONT_0(\Base)\) correspond to
\emph{continuous fields of Hilbert spaces} over~\(\Base\) in
the sense of~\cite{Dixmier:Cstar-algebras}.  The continuous
field corresponding to \(L^2(\Grd)\) has as its fibres the
Hilbert spaces \(L^2(\Grd_\base,\dd\lambda_\base)\) for
\(\base\in\Base\).  The topology on \(\bigsqcup_{\base\in\Base}
L^2(\Grd_\base,\dd\lambda_\base)\) is defined so that
\(\CONT_c(\Grd^{(1)})\) is dense in the resulting Banach space
of \(\CONT_0\)\nb-sections.

If \(\Grd=G\ltimes X\) for a locally compact group~\(G\) and a
locally compact \(G\)\nb-space~\(X\), then
\[
L^2(G\ltimes X)\cong L^2(G)\otimes \CONT_0(X)
\]
with~\(G\) acting diagonally.  Disregarding the group action,
the continuous field of Hilbert spaces corresponding to
\(L^2(G\ltimes X)\) is the trivial bundle with fibre \(L^2(G)\)
everywhere.  In contrast, for a general groupoid the continuous
field of Hilbert spaces \(L^2(\Grd)\) need not even be locally
trivial.

\subsection{The space of Fredholm operators}
\label{sec:Fred_space}

We can now define the \emph{set} of Fredholm operators we are
interested in.  We restrict attention to \emph{essentially
  unitary} operators, that is, operators whose image in the
Calkin algebra is unitary.

\begin{definition}
  \label{def:Fred_Hilbert}
  If~\(\Hils\) is a Hilbert space, then \(\Fred_\Hils\) or
  \(\Fred(\Hils)\) denotes the set of contractive, essentially
  unitary operators on~\(\Hils\):
  \[
  \Fred_\Hils \defeq \{T\in\Bound(\Hils)\mid
  \norm{T}\le 1,\quad
  TT^*-\ID_\Hils,\ T^*T-\ID_\Hils\in \Comp(\Hils)\}.
  \]
  If \(\Hils=(\Hils_x)_{x\in X}\) is a bundle of Hilbert spaces
  over some space~\(X\), we let
  \[
  \Fred_\Hils\defeq \bigsqcup_{x\in X} \Fred_{\Hils_x},
  \]
  equipped with the obvious map to~\(X\).  Here~\(\bigsqcup\)
  denotes the disjoint union.
\end{definition}

It is a subtle problem to topologise~\(\Fred_\Hils\).  It may
be impossible to define a good norm topology.  To understand
this, suppose first that our bundle of Hilbert spaces is trivial
with fibre~\(\Hils_0\).  Equivalently, we are dealing with the
Hilbert module \(\CONT_0(X,\Hils_0)\).  Then \(\Fred(\Hils) =
X\times \Fred(\Hils_0)\), and we may indeed equip
\(\Fred(\Hils_0)\) with the norm topology and then take the
product topology on \(\Fred(\Hils)\).  \emph{But this topology
  depends on the trivialisation.}

Two trivialisations differ by a unitary operator on
\(\CONT_0(X,\Hils_0)\).  Such unitary operators correspond to
\Star{}strongly continuous maps \(U\colon
X\to\Unitary(\Hils_0)\), where \(\Unitary(\Hils_0)\) denotes
the group of unitary operators on~\(\Hils_0\).  The
\Star{}strong topology means that both \(U(x)\xi\) and
\(U(x)^*\xi\) are continuous maps \(X\to\Hils_0\) for any
\(\xi\in\Hils_0\).  Unless~\(\Hils_0\) is finite-dimensional,
there are \Star{}strongly continuous maps
\(X\to\Unitary(\Hils_0)\) that are not norm-continuous.
Twisting our trivialisation by such a map yields a different
norm topology on~\(\Fred_\Hils\).

As a result, the norm topology on~\(\Fred_\Hils\) is only
well-defined if we have a preferred trivialisation to start
with.  While this is clearly the case for groupoids of the form
\(G\ltimes X\) for a group~\(G\), the bundle of Hilbert spaces
\(L^2(\Grd)\) need not even be locally trivial in general.

In order to construct a more canonical topology
on~\(\Fred_\Hils\), we first topologise bounded subsets in the
trivial bundles \(X\times \Bound(\Hils_0)\) and \(X\times
\Comp(\Hils_0)\) of bounded and compact operators.

On \(X\times \Bound(\Hils_0)\), we take the \Star{}strong
topology.  Thus a net \((x_\alpha,T_\alpha)\) in \(X\times
\Bound(\Hils_0)\) converges towards \((x,T)\) if and only if
\(\lim x_\alpha=x\) and \(\lim T_\alpha\xi=T\xi\) and \(\lim
T_\alpha^*\xi=\xi\) for all \(\xi\in\Hils_0\).  As long as
\((T_\alpha)\) is bounded, this topology is invariant under
left or right pointwise multiplication with \Star{}strongly
continuous maps \(X\to\Unitary(\Hils_0)\).  We have restricted
attention to Fredholm operators with \(\norm{T}\le1\) because
the multiplication is jointly continuous for the \Star{}strong
topology only on bounded subsets of \(\Bound(\Hils_0)\).

On \(X\times \Comp(\Hils_0)\), we take the norm topology, that
is, a net \((x_\alpha,T_\alpha)\) in \(X\times \Comp(\Hils_0)\)
converges towards \((x,T)\) if and only if \(\lim x_\alpha=x\)
and \(\lim {}\norm{T_\alpha-T}=0\).  The crucial point is that
on compact operators the norm topology is invariant under
multiplication by \Star{}strongly continuous maps \(U\colon
X\to\Unitary(\Hils_0)\).  For this, we must show that a net
\((x_\alpha,U(x_\alpha)\cdot T_\alpha)\) in \(X\times
\Comp(\Hils_0)\) converges to \((x,U(x)\cdot T)\) if
\((x_\alpha,T_\alpha)\) converges to \((x,T)\).  We estimate
\begin{align*}
  \norm{U(x_\alpha)T_\alpha- U(x)T}
  &\le \bigl\Vert U(x_\alpha)(T_\alpha - T) \bigr\rVert
  + \bigl\Vert \bigl(U(x_\alpha)- U(x)\bigr)T \bigr\rVert
  \\&\le \norm{T_\alpha-T}
  + \bigl\Vert \bigl(U(x_\alpha)- U(x)\bigr)T \bigr\rVert.
\end{align*}
The first term goes to~\(0\) by assumption.  To see that the
second one goes to~\(0\), we take \(\varepsilon>0\) and split
\(T=F+K\) with a finite-rank operator~\(F\) and
\(\norm{K}<\varepsilon\).  Then
\begin{align*}
  \bigl\Vert \bigl(U(x_\alpha)- U(x)\bigr)T \bigr\rVert
  &\le \bigl\Vert \bigl(U(x_\alpha)- U(x)\bigr)F \bigr\rVert
  + \bigl\Vert \bigl(U(x_\alpha)- U(x)\bigr)K \bigr\rVert
  \\&\le \bigl\Vert \bigl(U(x_\alpha)- U(x)\bigr)F \bigr\rVert
  + 2\varepsilon.
\end{align*}
The first term goes to~\(0\) because~\(F\) has finite rank and
\(U(x_\alpha)\to U(x)\) strongly.

As a result, the norm topology on \(X\times\Comp(\Hils_0)\) is
invariant under left multiplication by \Star{}strongly
continuous families of unitary operators.  A similar argument
deals with right multiplication.

So far, we have only considered operators on trivial Hilbert
modules.  Now let~\(\Hilm\) be any countably generated Hilbert
\(\CONT_0(X)\)-module.  By the (non-equivariant) Kasparov
Stabilisation Theorem, there exists an adjointable isometry
\[
V\colon \Hilm \xrightarrow{\textup{can}}
\Hilm\oplus \CONT_0(X,\ell^2\N) \cong
\CONT_0(X,\ell^2\N).
\]
View~\(\Hilm\) as a continuous field of Hilbert spaces
over~\(X\) with fibres~\(\Hilm_x\) for \(x\in X\).  Then~\(V\)
restricts to isometries on the fibres \(V(x)\colon
\Hilm_x\to\ell^2\N\).

\begin{definition}
  \label{def:topologise_operators}
  The \emph{\Star{}strong topology} on \(\bigsqcup
  \Bound(\Hilm_x)\) is defined as follows: a net
  \((x_\alpha,T_\alpha)\) converges in \(\bigsqcup
  \Bound(\Hilm_x)\) if and only if
  \((x_\alpha,V(x_\alpha)T_\alpha V(x_\alpha)^*)\) converges in
  \(X\times\Bound(\ell^2\N)\) in the \Star{}strong topology.

  Similarly, the \emph{norm topology} on \(\bigsqcup
  \Comp(\Hilm_x)\) is defined as follows: a net
  \((x_\alpha,T_\alpha)\) converges in \(\bigsqcup
  \Comp(\Hilm_x)\) if and only if
  \((x_\alpha,V(x_\alpha)T_\alpha V(x_\alpha)^*)\) converges in
  \(X\times\Comp(\ell^2\N)\) in the norm topology.
\end{definition}

If~\(W\) is another isometry as above, then there exists a
unitary operator~\(U\) on \(\CONT_0(X,\ell^2\N)\) with
\(W=UV\).  Since the \Star{}strong topology on the bounded
operators and the norm topology on the compact operators are
invariant under this change of trivialisation, the two
topologies in Definition~\ref{def:topologise_operators} are
well-defined.

\begin{lemma}
  The \Cstar{}algebra of compact operators on~\(\Hilm\) agrees
  with the \Cstar{}algebra of norm-continuous sections of the
  bundle \(\bigsqcup \Comp(\Hilm_x)\) that vanish at infinity.
\end{lemma}

\begin{proof}
  This is clear for trivial bundles of Hilbert spaces and
  follows in general from the Stabilisation Theorem.
\end{proof}

\begin{definition}
  \label{def:Fred_topology}
  Let \(\Hils=(\Hils_x)_{x\in X}\) be a continuous field of Hilbert spaces
  over~\(X\).  We equip the set \(\Fred_\Hils=\bigsqcup_{x\in X}
  \Fred_{\Hils_x}\) with the topology defined by the map
  \[
  \Fred_\Hils \to \bigsqcup_{x\in X} \Bound(\Hils_x)
  \times \bigsqcup_{x\in X} \Comp(\Hils_x)
  \times \bigsqcup_{x\in X} \Comp(\Hils_x),
  \qquad F\mapsto (F,1-FF^*,1-F^*F).
  \]
  That is, a net \((x_\alpha,F_\alpha)\) in~\(\Fred_{\Hils_x}\) converges to
  \((x,F)\) if and only if
  \begin{itemize}
  \item \(x_\alpha\) converges to~\(x\),

  \item \(F_\alpha\) converges to~\(F\) \Star{}strongly,

  \item \(1-F_\alpha F_\alpha^*\) converges to \(1-FF^*\) in norm,

  \item \(1-F_\alpha^* F_\alpha\) converges to \(1-F^*F\) in norm.
  \end{itemize}
  A similar definition applies to the space of Fredholm
  operators between the fibres of two different continuous
  fields of Hilbert spaces.
\end{definition}

Now we consider the continuous field of Hilbert spaces
\(L^2(\Grd)^\infty\) over the object space~\(\Base\) of a
locally compact groupoid~\(\Grd\) with Haar system.  We
abbreviate
\[
\Fred_\Grd \defeq \Fred\bigl(L^2(\Grd)^\infty\bigr).
\]
Since~\(\Grd\) acts continuously on \(L^2(\Grd)^\infty\), the
induced action of~\(\Grd\) on \(\Comp(L^2(\Grd)^\infty)\) by
conjugation is continuous for the norm topology, and the
induced action on the unit ball of \(\Bound(L^2(\Grd)^\infty)\)
is strongly continuous and hence \Star{}strongly continuous.
(Here we also use that the multiplication map is jointly
continuous for the strong topology on bounded subsets of
\(\Bound(L^2(\Grd)^\infty)\).)  Therefore, the induced action
of~\(\Grd\) on~\(\Fred_\Grd\) by conjugation on the fibres is
continuous in the topology described above.

\begin{definition}
  \label{def:support_Fred}
  Let \(f\colon X\to \Fred_\Hils\) be a
  \(\Grd\)\nb-equivariant continuous map.  Its \emph{support}
  is the closure of the set of \(x\in X\) where \(f(x)\) is not
  unitary.
\end{definition}

\begin{theorem}
  \label{the:RK_Fred_general}
  Let~\(\Grd\) be a locally compact groupoid with Haar system,
  form~\(\Fred_\Grd\) as above.  There are natural isomorphisms
  between
  \begin{itemize}
  \item \(\RK^0_\Grd(X)\) and the group of homotopy classes of
    \(\Grd\)\nb-equivariant continuous maps \(f\colon
    X\to\Fred_\Grd\);

  \item \(\K^0_\Grd(X)\) and the group of homotopy classes of
    \(\Grd\)\nb-equivariant continuous maps \(f\colon
    X\to\Fred_\Grd\) with \(\Grd\)\nb-compact support;

  \item \(\RK^0_{\Grd,Y}(X)\) and the group of homotopy classes
    of \(\Grd\)\nb-equivariant continuous maps \(f\colon
    X\to\nobreak\Fred_\Grd\) with \(Y\)\nb-compact support;
    here~\(Y\) is another \(\Grd\)\nb-space and \(p\colon X\to
    Y\) is a continuous \(\Grd\)\nb-equivariant map.

  \end{itemize}
\end{theorem}

The description of \(\K^0_\Grd(X)\) in
Theorem~\ref{the:RK_Fred_general} is a special case of
\cite{Tu-Xu-Laurent-Gengoux:Twisted_K}*{Theorem 3.14}, which
describes equivariant \emph{twisted} \(\K\)\nb-theory in terms
of Fredholm operators.

\subsection{Proof of the theorem}
\label{sec:proof_RK_Fred}

Let \(p\colon X\to Y\) be a \(\Grd\)\nb-space over~\(Y\).  We
are going to prove the assertion about \(\RK^0_{\Grd,Y}(X)\) in
Theorem~\ref{the:RK_Fred_general}.  We first recall Kasparov's
definition of
\[
\RK^0_{\Grd,Y}(X) \defeq
\KK^{\Grd\ltimes Y}_0\bigl(\CONT_0(Y),\CONT_0(X)\bigr).
\]
A cycle for \(\RK^0_{\Grd,Y}(X)\) consists of a
\(\Z/2\)\nb-graded \(\Grd\)\nb-equivariant Hilbert
module~\(\Hilm\) over~\(\CONT_0(X)\), a \(\Grd\ltimes
Y\)\nb-equivariant \Star{}homomorphism~\(\varphi\) from
\(\CONT_0(Y)\) to the \Cstar{}algebra of adjointable, grading
preserving operators on~\(\Hilm\), and a self-adjoint, odd,
almost \(\Grd\)\nb-equivariant operator~\(F\) on~\(\Hilm\),
such that \((F^2-1)\varphi(h)\) and \([F,\varphi(h)]\) are
compact for all \(h\in\CONT_0(Y)\).  Two cycles give the same
element in \(\RK^0_{\Grd,Y}(X)\) if and only if they are
homotopy equivalent, where homotopies are, by definition,
cycles for \(\RK^0_{\Grd,Y}(X\times[0,1])\).

The \(\Grd\ltimes Y\)-linearity completely
determines~\(\varphi\): it must be given by pointwise
multiplication: \(\varphi(f)(\xi)=\xi\cdot (f\circ p)\) for all
\(\xi\in\Hilm\), \(f\in\CONT_0(Y)\).  Hence we can omit this
part of the data and are left with a pair \((\Hilm,F)\).  Since
\([F,\varphi(h)]=0\) for all \(h\in\CONT_0(Y)\), we can also
forget about one of the compactness conditions.  Since~\(\Grd\)
acts properly on~\(X\), any almost \(\Grd\)\nb-equivariant
operator on~\(\Hilm\) has a compact perturbation that is
exactly \(\Grd\)\nb-equivariant (see
\cite{Tu:Novikov}*{\S6.3}).  Therefore, we may restrict
attention to cycles with \(\Grd\)\nb-equivariant~\(F\).  In
addition, we can achieve \(\norm{F}\le 1\) by functional
calculus.  We assume that this is the case from now on.

Let~\(\Hilm^\pm\) be the even and odd subspaces of~\(\Hilm\)
with respect to the \(\Z/2\)\nb-grading.  The block matrix
decomposition of~\(F\) is
\[
F= \begin{pmatrix}0&U^*\\U&0\end{pmatrix}
\]
for an adjointable operator \(U\colon \Hilm^+\to\Hilm^-\) with
\(\norm{U}\le 1\) because~\(F\) is odd, self-adjoint and
contractive.  Hence we may replace \((\Hilm,F)\) by
\((\Hilm^+,\Hilm^-,U)\), where~\(\Hilm^\pm\) are
\(\Grd\)\nb-equivariant Hilbert modules over~\(\CONT_0(X)\) and
\(U\colon \Hilm^+\to\Hilm^-\) is a \(\Grd\)\nb-equivariant
adjointable operator of norm~\(1\), such that \((1-U^*U)\cdot
\varphi(f\circ p)\) and \((1-UU^*)\cdot \varphi(f\circ p)\) are
compact operators on \(\Hilm^+\) and~\(\Hilm^-\), respectively,
for all \(f\in\CONT_0(Y)\).

The Equivariant Stabilisation Theorem~\ref{the:stabilisation}
asserts that there are \(\Grd\)\nb-equivariant unitary
operators
\[
\Hilm^\pm\oplus L^2(\Grd\ltimes X)^\infty \cong L^2(\Grd\ltimes X)^\infty.
\]
Since addition of the degenerate cycle \(\ID\colon
L^2(\Grd\ltimes X)^\infty\to L^2(\Grd\ltimes X)^\infty\) leads
to a homotopic cycle, we may restrict attention to cycles whose
underlying Hilbert module is \(L^2(\Grd\ltimes X)^\infty\).
Thus only the operator~\(U\) remains as data.

\begin{lemma}
  \label{lem:homotopy_unstable}
  Two cycles \(U_0,U_1\colon L^2(\Grd\ltimes X)^\infty\to
  L^2(\Grd\ltimes X)^\infty\) for \(\RK^0_{\Grd,Y}(X)\) have
  the same class in \(\RK^0_{\Grd,Y}(X)\) if and only if there
  is a cycle~\(U\) on \(L^2(\Grd\ltimes X\times[0,1])^\infty\)
  that restricts to \(U_0\) and~\(U_1\) at the endpoints.
\end{lemma}

\begin{proof}
  Let \(\Hils\defeq L^2(\Grd\ltimes X)\otimes L^2([0,1])\),
  where we use the Lebesgue measure on \([0,1]\).  This
  \(\Grd\)\nb-equivariant Hilbert module over \(\CONT_0(X)\) is
  unitarily equivalent to \(L^2(\Grd\ltimes X)^\infty\), so
  that we may replace \(U_0\) and~\(U_1\) by operators
  on~\(\Hils\).  By the definition of our equivalence relation,
  the cycles \(U_0\) and~\(U_1\) are homotopic via a homotopy
  that is realised on some \(\Grd\)\nb-equivariant Hilbert
  module~\(\Hilm\) over \(\CONT_0(X\times[0,1])\).  Using the
  Equivariant Stabilisation Theorem, we get a homotopy between
  \(U_0\oplus\ID\) and \(U_1\oplus\ID\) realised on the Hilbert
  module \(\Hilm\oplus \Hils\otimes \CONT([0,1]) \cong
  \Hils\otimes \CONT([0,1])\).

  We identify \(\Hils\cong \Hils\oplus\Hils\) using the unitary
  operator
  \[
  V\colon L^2([0,1]) =
  L^2([0,\nicefrac{1}{2}]) \oplus
  L^2([\nicefrac{1}{2},1]) \xrightarrow[\cong]{S_1\oplus S_2}
  L^2([0,1]) \oplus L^2([0,1]),
  \]
  where
  \begin{alignat*}{2}
    S_1(f)(t) &\defeq \sqrt{2} f(t/2) &\qquad&
    \text{for \(f\in L^2([0,\nicefrac{1}{2}])\),}\\
    S_2(f)(t) &\defeq \sqrt{2} f\bigl((t+1)/2\bigr) &\qquad&
    \text{for \(f\in L^2([\nicefrac{1}{2},1])\).}
  \end{alignat*}
  Conjugating by~\(V\), the homotopy above yields one between
  \(V^*(U_0\oplus\ID_\Hils)V\) and \(V^*(U_1\oplus\ID_\Hils)V\)
  that is realised on \(\Hils\otimes\CONT([0,1])\).  It
  remains, therefore, to find a homotopy between
  \(V^*(U_t\oplus\ID_\Hils)V\) and~\(U_t\).  For this, we use
  the \Star{}strongly continuous family of unitary operators
  \[
  V_s\colon L^2([0,1]) = L^2([0,s])\oplus L^2([s,1]) \cong
  L^2([0,1])\oplus L^2([0,1])
  \]
  for \(s\in [1/2,1]\) defined like \(S_1\oplus S_2\).  The
  operators \(V_s^*(U_t\oplus\ID_\Hils)V_s\) for \(t\in[0,1]\)
  provide the desired homotopy between \(U_t\oplus\ID_\Hils\)
  and~\(U_t\).
\end{proof}

The homotopies in the proof of
Lemma~\ref{lem:homotopy_unstable} are only \Star{}strongly
continuous.  Such homotopies are much easier to accomplish than
norm-continuous homotopies.

Now we identify \(L^2(\Grd\ltimes X)^\infty\) with the space of
\(\CONT_0\)\nb-sections of a continuous field of Hilbert spaces
over~\(X\) and replace \(U\colon L^2(\Grd\ltimes X)^\infty\to
L^2(\Grd\ltimes X)^\infty\) by a family \((U_x)_{x\in X}\) of
operators between the fibres of this continuous field.

Since \((1-UU^*)\varphi(f\circ p)\) and
\((1-U^*U)\varphi(f\circ p)\) are compact for all
\(f\in\CONT_0(Y)\), the operators \(1-U_xU_x^*\) and
\(1-U_x^*U_x\) are compact for all \(x\in X\).  We also have
\(\norm{U_x}\le1\) for all \(x\in X\).  Thus \(x\mapsto U_x\)
is a map from~\(X\) to \(\Fred_\Grd\).  The equivariance
of~\(U\) is equivalent to the \(\Grd\)\nb-equivariance of
this map.

Since \((U_x)_{x\in X}\) defines an adjointable operator on
\(L^2(\Grd\ltimes X)^\infty\), the map \(x\mapsto U_x\) is
\Star{}strongly continuous.  Since \((1-UU^*)\varphi(f\circ
p)\) and \((1-U^*U)\varphi(f\circ p)\) are compact for all
\(f\in\CONT_0(Y)\) and continuity is a local issue, the maps
\(x\mapsto 1-U_xU_x^*\) and \(x\mapsto 1-U_x^*U_x\) are
norm-continuous.  Thus the map \(x\mapsto U_x\) is continuous
for the topology in Defition~\ref{def:Fred_topology}.  Finally,
the compactness of \((1-UU^*)\varphi(f\circ p)\) and
\((1-U^*U)\varphi(f\circ p)\) implies that
\(\norm{(1-U_xU_x^*)}\cdot f\bigl(p(x)\bigr)\) and
\(\norm{(1-U_x^*U_x)}\cdot f\bigl(p(x)\bigr)\) are
\(\CONT_0\)\nb-functions on~\(X\).  Equivalently, the subsets
\[
\bigl\{x\in X \bigm|
\text{\(\norm{1-U_x^*U_x}\ge\varepsilon\) or
  \(\norm{1-U_xU_x^*} \ge\varepsilon\)} \bigr\}
\]
are \(Y\)\nb-compact for all \(\varepsilon>0\).  This is not
quite the support condition that we want, but we can improve
this condition using functional calculus.

Let \(f\colon \R\to[0,1]\) be a continuous function with
\(f(t)=t^{-1/2}\) for \(t\ge\nicefrac{1}{2}\) and let
\(U'\defeq U\cdot f(U^*U)\), so that \(U'_x= U_x\cdot
f(U_x^*U_x)\) for all \(x\in X\).  This yields a homotopic
cycle for \(\RK^*_{\Grd,Y}(X)\) via a linear homotopy between
\(U\) and~\(U'\) because~\(U\) is essentially unitary and
\(f(1)=1\).  The operator~\(U'_x\) is unitary where both
\(\norm{1-U_x^*U_x}\le\nicefrac{1}{2}\) and
\(\norm{1-U_xU_x^*}\le\nicefrac{1}{2}\).  Hence \((U'_x)_{x\in
  X}\) has \(Y\)\nb-compact support.

Thus a cycle for \(\RK^0_{\Grd,Y}(X)\) yields a
\(\Grd\)\nb-equivariant continuous map \({X\to\Fred_\Grd}\)
with \(Y\)\nb-compact support.  Since we can apply the same
construction to homotopies, Lemma~\ref{lem:homotopy_unstable}
yields that homotopic cycles for \(\RK^0_{\Grd,Y}(X)\) yield
homotopic maps \({X\to\Fred_\Grd}\).  Conversely, a
\(\Grd\)\nb-equivariant continuous map
\(X\to\Fred_{\Grd\ltimes X}\) with \(Y\)\nb-compact support
yields a cycle for \(\RK^0_{\Grd,Y}(X)\), and homotopic maps
yield homotopic cycles.  Thus \(\RK^0_{\Grd,Y}(X)\) agrees with
the set of homotopy classes of \(\Grd\)\nb-equivariant
continuous maps \(X\to\Fred_\Grd\) with \(Y\)\nb-compact
support as asserted in Theorem~\ref{the:RK_Fred_general}.

The description of \(\RK^0_{\Grd,Y}(X)\) contains
\(\RK^0_\Grd(X)=\RK^0_{\Grd,X}(X)\) as a special case; notice
that if~\(p\) is the identity map on~\(X\), then any closed
subset of~\(X\) is \(X\)\nb-compact.  Thus we get the assertion
about \(\RK^0_\Grd(X)\) in Theorem~\ref{the:RK_Fred_general}.

Now we turn to the group \(\K^0_\Grd(X)\), which we describe as
\(\KK_0\bigl(\C,\Grd\ltimes \CONT_0(X)\bigr)\).  Cycles for
\(\KK_0\bigl(\C,\Grd\ltimes \CONT_0(X)\bigr)\) consist of two
Hilbert modules~\(\tilde\Hilm^\pm\) over \(\Grd\ltimes
\CONT_0(X)\) together with an adjointable operator
\(\tilde{U}\colon \Hilm^+\to\Hilm^-\) such that
\(1-\tilde{U}\tilde{U}^*\) and \(1-\tilde{U}^*\tilde{U}\) are
compact.  The \Cstar{}category of Hilbert modules over
\(\Grd\ltimes \CONT_0(X)\) is equivalent to the
\Cstar{}category of \(\Grd\)\nb-equivariant Hilbert modules
over \(\CONT_0(X)\); this is implicit in
\cite{Tu:Novikov}*{Proposition 6.24}.  Hence we may replace
\(\tilde{\Hilm}^\pm\) by \(\Grd\)\nb-equivariant Hilbert
modules~\(\Hilm^\pm\) over \(\CONT_0(X)\) and~\(\tilde{U}\) by
a \(\Grd\)\nb-equivariant adjointable operator \(U\colon
\Hilm^+\to\Hilm^-\).  Applying the Equivariant Stabilisation
Theorem~\ref{the:stabilisation}, we reduce to
\(\Hilm^\pm=L^2(\Grd\ltimes X)^\infty\).

A \(\Grd\)\nb-equivariant adjointable operator~\(T\) on
\(L^2(\Grd\ltimes X)^\infty\) corresponds to a compact operator
on \(C^*(\Grd\ltimes X)^\infty\) if and only if \(T\cdot
\varphi(c)\) is compact, where~\(c\) is a cut-off function.
This is a variant of \cite{Tu:Novikov}*{Proposition 6.24}.  Now
the argument is essentially the same as above.

\subsection{Comparing the support conditions}
\label{sec:Fred_compare_support}

The descriptions of
\[
\RK^0_\Grd(X),
\qquad
\RK^0_{\Grd,Y}(X),
\quad\text{and}\quad
\K^0_\Grd(X)
\]
in Theorem~\ref{the:RK_Fred_general} differ only in the support
conditions.  This justifies calling \(\RK^0_{\Grd,Y}(X)\) the
\(\K\)\nb-theory of~\(X\) with \(Y\)\nb-compact and
\(\K^0_\Grd(X)\) the \(\K\)\nb-theory of~\(X\) with
\(\Grd\)\nb-compact support, whereas \(\RK^0_\Grd(X)\) is the
\(\K\)\nb-theory of~\(X\) without support restriction.

\begin{proposition}
  \label{pro:Y-compact_proper_change_of_Y}
  There are canonical maps
  \[
  \K^*_\Grd(X) \to \RK^*_{\Grd,Y}(X) \to \RK^*_\Grd(X).
  \]
  The first is an isomorphism if~\(Y\) is
  \(\Grd\)\nb-compact.  The second map is an isomorphism once
  the map \(X\to Y\) is proper.  Both maps are isomorphisms
  if~\(X\) is \(\Grd\)\nb-compact.

  A continuous \(\Grd\)\nb-map \(h\colon Y\to Y'\) induces a
  map \(\RK^*_{\Grd,Y'}(X) \to \RK^*_{\Grd,Y}(X)\), which is an
  isomorphism if~\(h\) is proper.
\end{proposition}

\begin{proof}
  It suffices to treat the even \(\K\)\nb-groups, where we
  can use Theorem~\ref{the:RK_Fred_general}.  Thus it suffices
  to compare the various support conditions.

  We claim that any \(\Grd\)\nb-compact subset of~\(X\) is
  \(Y\)\nb-compact.  Let \(K\subseteq X\) and \(L\subseteq Y\)
  be compact and \(g\in\Grd^{(1)}\), then \(g\cdot K\cap
  p^{-1}(L)\neq\emptyset\) if and only if \(g\cdot p(K)\cap
  L\neq\emptyset\); since~\(Y\) is proper, the set of
  \(g\in\Grd^{(1)}\) for which this happens is compact, so that
  \(\Grd\cdot K\cap p^{-1}(L)\) is compact.  This means that
  \(p\) restricts to a proper map on \(\Grd\cdot K\), so that
  all \(\Grd\)\nb-invariant closed subsets of \(\Grd\cdot K\)
  are \(Y\)\nb-compact.

  Since \(\Grd\)\nb-compact subsets of~\(X\) are
  \(Y\)\nb-compact, we get a map \(\K^0_\Grd(X) \to
  \RK^0_{\Grd,Y}(X)\).  If~\(Y\) is \(\Grd\)\nb-compact, say,
  \(Y=\Grd\cdot L\), then any \(Y\)\nb-compact
  \(\Grd\)\nb-invariant subset~\(A\) of~\(X\) is
  \(\Grd\)\nb-compact because it is contained in \(\Grd\cdot
  \bigl(p^{-1}(L)\cap A\bigr)\).  Therefore, \(\K^0_\Grd(X)
  \cong \RK^0_{\Grd,Y}(X)\) if~\(Y\) is \(\Grd\)\nb-compact.
  If~\(X\) is \(\Grd\)\nb-compact, the same argument shows that
  all our support conditions are vacuously satisfied, so that
  all three groups coincide by
  Theorem~\ref{the:RK_Fred_general}.
  Theorem~\ref{the:RK_Fred_general} yields a canonical map
  \(\RK^0_{\Grd,Y}(X) \to \RK^0_\Grd(X)\), forgetting the
  support restriction.

  Finally, consider a map \(h\colon Y\to Y'\).  Then any
  \(Y'\)\nb-compact subset of~\(X\) is also \(Y\)\nb-compact,
  and the converse holds if and only if~\(h\) is proper.
  Therefore, we get a map \(\RK^0_{\Grd,Y'}(X) \to
  \RK^0_{\Grd,Y}(X)\), and it is an isomorphism if~\(h\) is
  proper.
\end{proof}

\subsection{Changing the topology}
\label{sec:topology}

Now we restrict attention to \emph{group} actions on spaces.
In this case, we do have a canonical trivialisation
\[
L^2(G\ltimes X)^\infty\cong \CONT_0(X,L^2(G)^\infty),
\]
with~\(G\) acting diagonally.  Hence the norm topology on the
bundle of Fredholm operators is meaningful.  In this section,
we let~\(\Fred_G\) be the space of essentially unitary
operators on \(\Hils_G\defeq L^2(G)^\infty\), equipped with the
\emph{norm topology}.

\begin{theorem}
  \label{the:RK_Fred_group}
  Let~\(G\) be a locally compact group.  There are natural
  isomorphisms between
  \begin{itemize}
  \item \(\RK^0_G(X)\) and the group of homotopy classes of
    \(G\)\nb-equivariant \textup{(}norm\textup{)} continuous
    maps \(f\colon X\to\Fred_G\);

  \item \(\K^0_G(X)\) and the group of homotopy classes of
    \(G\)\nb-equivariant continuous maps \(f\colon
    X\to\Fred_G\) with \(G\)\nb-compact support;

  \item \(\RK^0_{G,Y}(X)\) and the group of homotopy classes of
    \(G\)\nb-equivariant continuous maps \(f\colon
    X\to\Fred_G\) with \(Y\)\nb-compact support; here~\(Y\)
    is a proper \(\Grd\)\nb-space and \(p\colon X\to Y\) is a
    \(G\)\nb-equivariant continuous map.

  \end{itemize}
\end{theorem}

\begin{corollary}
  \label{cor:unitary_contractible}
  The space of norm-continuous \(G\)\nb-equivariant maps
  from~\(X\) to the unitary group \(\Unitary(\Hils_G)\) is
  contractible.
\end{corollary}

\begin{proof}
  Since maps \(X\to\Fred_G\) that are unitary everywhere
  represent zero in \(\RK^0_G(X)\),
  Theorem~\ref{the:RK_Fred_group} shows that this space of maps
  is connected.  Since we may replace~\(X\) by
  \(X\times\Sphere^n\) for \(n\in\N\), this space is weakly
  contractible as well.  We may replace unitaries by
  invertibles without changing the homotopy type.  This yields
  an open subset of the \Cstar{}algebra of \(G\)\nb-equivariant
  norm-continuous maps \(X\to\Bound(\Hils_G)\).  Open subsets
  of Banach spaces have the homotopy type of a CW-complex by
  \cite{Lundell-Weingram:Topology_CW}*{Corollary 5.5 in
    Chapter~4}.  Hence they are contractible once they are
  weakly contractible.
\end{proof}

The only difference between Theorems \ref{the:RK_Fred_general}
and~\ref{the:RK_Fred_group} is the topology on the bundle of
Fredholm operators.

The proof of Theorem~\ref{the:RK_Fred_group} will occupy the
remainder of this section.  We concentrate on the case
\(\RK^0_{G,Y}(X)\) because it contains \(\RK^0_G(X)\) as a
special case and \(\K^0_G(X)\) is analogous.

The first part of the argument is similar to the proof of
\cite{Emerson-Meyer:Descent}*{Proposition 13}.  Let~\(D\ns\) be
the \Cstar{}algebras of \(G\)\nb-equivariant maps
\(X\to\Bound(\Hils_G)\) that are \Star{}strongly continuous,
and let \(D\nc\subseteq D\ns\) be the subalgebra of
norm-continuous maps.

\begin{lemma}
  \label{lem:Eilenberg_swindle}
  \(\K_*(D\nc)=0\) and \(\K_*(D\ns)=0\).
\end{lemma}

\begin{proof}
  Both assertions are proved using an Eilenberg swindle.  We
  carry this out only for \(\K_0(D\nc)\) because the other
  three cases \(\K_1(D\nc)\), \(\K_0(D\ns)\), and
  \(\K_1(D\ns)\) are analogous.  Since \(D\nc\) is unital and
  matrix stable, it suffices to study idempotents in~\(D\nc\).
  Let~\(p\) be such an idempotent and write \(p^\bot\defeq
  1-p\).  Choose an isomorphism \(\Hils_G\cong\Hils_G^\infty\)
  and use this to transform the diagonal operator \(\infty\cdot
  p \defeq p\oplus p\oplus p\oplus p\oplus \dotsb\) into an
  operator on~\(\Hils_G\).  Notice that this is again a
  norm-continuous function \(X\to\Bound(\Hils_G)\) if~\(p\) is.
  The projections \(p\oplus \infty\cdot p\) and
  \(0\oplus\infty\cdot p\) are equivalent via the unilateral
  shift operator on~\(\Hils_G^\infty\), which we view as a
  constant function \(X\to\Bound(\Hils_G^\infty)\).  Hence
  \([p]=0\) in \(\K_0(D\nc)\).
\end{proof}

Let~\(J_0\) be the algebra of norm-continuous
\(G\)\nb-equivariant maps \(X\to\Comp(\Hils_G)\) whose support
is \(Y\)\nb-compact.  This is a (non-closed) \Star{}ideal in
both \(D\nc\) and \(D\ns\).  Its closure is the
\Cstar{}algebra~\(J\) of all norm-continuous
\(\Grd\)\nb-equivariant maps \(f\colon X\to\Comp(\Hils_G)\) for
which the subsets \(\{x\in X\mid \norm{f(x)}>\varepsilon\}\)
are \(Y\)\nb-compact for all \(\varepsilon>0\).

Let \(Q\nc\) and \(Q\ns\) be the quotients of
\(D\nc\) and \(D\ns\) by the ideal~\(J\), respectively.

\begin{lemma}
  \label{lem:RKK_as_Kone}
  \(\K_{*+1}(Q\nc)\cong \K_{*+1}(Q\ns)\cong
  \RK^*_{G,Y}(X)\).
\end{lemma}

\begin{proof}
  The embedding \(D\nc\to D\ns\) induces a morphism of
  \Cstar{}algebra extensions
  \[
  \xymatrix{
    J \ar[r] \ar@{=}[d] & D\nc \ar[d] \ar[r] &
    Q\nc \ar[d] \\
    J \ar[r] & D\ns \ar[r] & Q\ns.
  }
  \]
  The vertical maps \(J\to J\) and \(D\nc\to D\ns\) induce
  isomorphisms on \(\K\)\nb-theory by
  Lemma~\ref{lem:Eilenberg_swindle}.  Now the Five Lemma
  applied to the long exact \(\K\)\nb-theory sequence yields
  \(\K_*(Q\nc)\cong \K_*(Q\ns)\).

  For the second isomorphism, we only prove \(\K_1(Q\ns)\cong
  \RK^0_{G,Y}(X)\), the other parity can be reduced to this by
  suspension.  The \Cstar{}algebra~\(Q\ns\) is unital and
  matrix stable, that is, \(Q\ns \cong \Mat_n(Q\ns)\).  Hence
  any cycle for \(\K_1(Q\ns)\) is represented by a unitary
  element in \(Q\ns\).  A lifting for such a unitary element
  to~\(D\ns\) is nothing but an adjointable operator~\(U\) on
  the Hilbert module \(\CONT_0(X,\Hils_G)\) with
  \[
  1-UU^*\in J \qquad\text{and}\qquad 1-U^*U\in J.
  \]
  These are exactly the cycles that describe \(\RK^0_{G,Y}(X)\).

  The equivalence relation on \(\K_1(Q\ns)\) is generated by
  stabilisation --~replacing \(U\) by \(V(U\oplus
  \ID_{\Hils_G})V^*\) for a unitary \(V\colon
  \Hils_G^2\to\Hils_G\)~-- and norm-continuous homotopy
  in~\(U\).  In Kasparov theory, these two relations are called
  addition of degenerate cycles and operator homotopy, and are
  known to generate the same equivalence relation as homotopy
  (in the non-equivariant case, this is
  \cite{Blackadar:K-theory}*{Theorem 18.5.3}; the same argument
  works for equivariant Kasparov theory).  Hence
  \(\K_1(Q\ns)\cong \RK^0_{G,Y}(X)\).
\end{proof}

The group \(\K_1(Q\nc)\) can be described similarly.  The only
difference to \(\K_1(Q\ns)\) is that its cycles are \emph{norm}
continuous \(G\)\nb-equivariant maps \(U\colon X\to\Fred_G\).
As in the proof of Theorem~\ref{the:RK_Fred_general}, we can
restrict attention to maps~\(f\) with \(Y\)\nb-compact support.
The equivalence relation on \(\K_1(Q\ns)\) is generated by two
moves: norm-continuous homotopy and stabilisation.  It remains
to prove, therefore, that for any norm-continuous
\(G\)\nb-equivariant map \(U\colon X\to\Fred_G\) with
\(Y\)\nb-compact support, there is a norm-continuous homotopy
between \(U\) and \(V(U\oplus \ID_{\Hils_G})V^*\) for the
standard unitary \(V\colon \Hils_G^2\to\Hils_G\).

As in more classical situations
(\cite{Jaenich:Vektorraumbuendel}), this follows from a
suitable generalisation of Kuiper's Theorem (see
Corollary~\ref{cor:unitary_contractible}).  In our case, we
need to know that the unitary group of~\(D\nc\) is
contractible.  The most general form of Kuiper's Theorem in the
literature is due to Joachim Cuntz and Nigel Higson
(\cite{Cuntz-Higson:Kuiper}, see also
\cite{Wegge-Olsen:K-theory}*{\S16} for a more detailed proof);
it asserts that the unitary group of a stable multiplier
algebra \(\Mult(A\otimes\Comp)\) is norm-contractible for any
\(\sigma\)\nb-unital \Cstar{}algebra~\(A\).  Unfortunately,
this does not cover the algebra~\(D\nc\) we need: only~\(D\ns\)
is a stable multiplier algebra.  Nevertheless, inspection shows
that the argument in~\cite{Cuntz-Higson:Kuiper} can be carried
over to~\(D\nc\).  We concentrate on those steps in the
argument that require change.  We modify the structure of the
argument slightly to simplify the application to Fredholm
operators, which is not considered in~\cite{Cuntz-Higson:Kuiper}.

We fix an isometry \(V = (V_1\ V_2)\colon \Hils_G^2\to\Hils_G\) as
above, that is, \(V_1\) and~\(V_2\) are isometries with
\(V_1V_1^*+V_2V_2^*=1\), and let \(P\defeq V_2V_2^*\).  We view
\(P,V_1,V_2\) as constant functions on~\(X\) and thus as elements of
\(D\nc\).

\begin{lemma}
  \label{lem:technical_isometry}
  Let \(\bar{U}\in Q\nc\) be unitary.  There are a lifting~\(U\!\)
  of~\(\bar{U}\!\) and isometries \(S\) and~\(T\) in \(D\nc\) such that:
  \begin{itemize}
  \item \(US\) is an isometry, and \(US\) and \(U(1-SS^*)\)
    have orthogonal range;

  \item \(V_2^*S=0\) and \(V_2^*T=0\), that is, the range projections of
    \(S\) and~\(T\) are orthogonal to~\(P\);

  \item \(\norm{T^*US}<1\).
  \end{itemize}
\end{lemma}

Before we prove this technical lemma, we use it to finish the proof of
Theorem~\ref{the:RK_Fred_group}.  As in~\cite{Cuntz-Higson:Kuiper}, we
call two projections \(P\) and~\(Q\) in a unital \Cstar{}algebra~\(A\)
\emph{strongly equivalent} and write \(P\approx Q\) if there is a
continuous path of unitaries \(U\in\CONT([0,1],A)\) with \(U_0=1_A\) and
\(U_1PU_1^*=Q\).  Two projections are strongly equivalent if and only if
they are homotopic (\cite{Blackadar:K-theory}*{Proposition 4.3.3}).
\cite{Cuntz-Higson:Kuiper}*{Lemma 1} asserts that two projections \(P\)
and~\(Q\) with \(\norm{PQ}<1\) are strongly equivalent once they are
equivalent.

Let \(R_S\), \(R_{US}\), and~\(R_T\) be the range projections
of the isometries \(S\), \(US\), and~\(T\), respectively.
Recall that~\(P\) is the range projection of~\(V_2\).  As range
projections of isometries, these projections are all equivalent
to~\(1\).  Since \(PR_S=PR_T=0\) and \(\norm{R_T R_{US}}<1\) by
Lemma~\ref{lem:technical_isometry},
\cite{Cuntz-Higson:Kuiper}*{Lemma 1} yields strong equivalences
\[
R_S\approx P,
\qquad
R_{US}\approx R_T \approx P.
\]
Thus we get continuous paths of unitaries
\(W_1\in\CONT([0,1],D\nc)\) and
\(W_2\in\CONT([0,1],D\nc)\) with \(W_1(0)=W_2(0)=1\) and
\[
W_1(1)R_{US}W_1(1)^* = P,
\qquad
W_2(1)R_SW_2(1)^* = P.
\]
This provides a homotopy between~\(U\) and \(U'\defeq
W_1(1)UW_2(1)^*\).  By construction,
\begin{multline*}
  U'P
  = W_1(1) UW_2(1)^*P
  = W_1(1) U R_S W_2(1)^*
  \\= W_1(1) R_{US} U W_2(1)^*
  = P W_1(1) U W_2(1)^*
  = P U'.
\end{multline*}
Hence
\[
U' = V \begin{pmatrix}U_1&0\\0&U_2\end{pmatrix}V^*
\]
for the unitary \(V\colon \Hils_G^2\to\Hils_G\) chosen above,
with essentially unitary \(U_1\) and~\(U_2\).  The usual
rotation homotopy that proves the commutativity of~\(\K_1\)
shows that \(V(U_1\oplus U_2)V^*\) is homotopic to
\(V(U_1U_2\oplus 1)V^*\).  But if two unitaries of the latter
form are stably homotopic, then they are already homotopic
because we can subsume any additional stabilisation in the
summand~\(1\).

Hence \(\K_1(Q\nc)\) is the set of homotopy classes of
unitaries in~\(Q\nc\), without need to stabilise.  Thus
\(\RK^0_G(X) \cong \K_1(Q\nc) \cong \K_1(Q\ns) \cong
[X,\Fred_G]\), as desired.

To finish the proof of Theorem~\ref{the:RK_Fred_group}, it
remains to prove Lemma~\ref{lem:technical_isometry}.

Splitting \(\ell^2\N\cong \ell^2\N\otimes\ell^2\N\), we can
embed \(\Bound(\ell^2\N) \otimes D\nc\) into \(D\nc\).  We
choose an infinite orthonormal sequence of rank-one projections
\((e_n)_{n\in\N}\) on \(\ell^2\N\) such that \(\sum e_n \otimes
1 = V_1V_1^*\).  Next we choose an approximate unit for the
ideal \(J\subseteq D\nc\) with \(u_{n+1}u_n = u_n\) for all
\(n\in\N\) that is \emph{locally uniformly continuous} as a set
of functions on~\(X\).  We omit the proof that such an
approximate unit exists.  We view \(e_n\otimes u_n\) as
elements of \(D\nc\).  By construction, \((e_n\otimes u_n)P=0\)
for all \(n\in\N\).  The local uniform continuity of \((u_n)\)
is needed for the following analogue of
\cite{Cuntz-Higson:Kuiper}*{Lemma 3}:

\begin{lemma}
  \label{lem:small_projections}
  Let \((u_n)_{n\in\N}\) and \((e_n)_{n\in\N}\) be as above,
  and let \(k\colon \N\to\N\) be an increasing function with
  \(\lim k(n)=\infty\), so that \((u_{k(n)})\) and
  \((e_{k(n)})\) are subsequences of \((u_n)\) and \((e_n)\).
  Let \(a\defeq \sum_{n\in\N} e_{k(n)}\otimes u_{k(n)}\).  Then
  there is an isometry \(T\in D\nc\) with \(aT = T\) and hence
  \(aTT^*=TT^*\) and \(TT^*\le a\).
\end{lemma}

\begin{proof}
  The proof is identical to the proof of
  \cite{Cuntz-Higson:Kuiper}*{Lemma 3}.  The assumption
  on~\(u_n\) implies that the sequence \((d_n)\) defined by
  \(d_n \defeq (u_{k(n-1)}-u_{k(n-2)})^{\nicefrac{1}{2}}\) is
  locally uniformly continuous, so that the sum \(T\defeq
  \sum_{n\in\N}^\infty v_n\otimes d_n\) used
  in~\cite{Cuntz-Higson:Kuiper} is norm-continuous; here
  \((v_n)\) is a sequence of partial isometries with carefully
  selected range and source projections.
\end{proof}

\begin{proof}[of Lemma~\ref{lem:technical_isometry}]
  Choose any lifting~\(U_0\) of~\(\bar{U}\).  First we
  modify~\(U_0\) and find an isometry~\(S_0\) to get the first
  property in the lemma.  Let \((e_n)\) and \((u_n)\) be as
  above, and let \(g_n\defeq e_n\otimes u_n\), viewed as
  elements of~\(D\nc\).  The sequence \((g_n)\) converges
  strictly to~\(0\) because \((e_n)\to0\) in the strict
  topology on \(\Bound(\ell^2\N)\).  Since \(1-U_0^*U_0\in J\),
  we can find a subsequence \(k(n)\) with
  \[
  \norm{(1-U_0^*U_0)g_{k(n)}} < \varepsilon\,2^{-n-1}.
  \]
  Lemma~\ref{lem:small_projections} yields an isometry~\(S_0\)
  with \(\sum g_{k(n)} \cdot S_0 = S_0\).  Since \(g_n P =
  (e_n\otimes u_n)P=0\) for all \(n\in\N\), the range
  of~\(S_0\) is orthogonal to~\(P\) as desired.  We have
  \[
  \norm{(1-U_0^*U_0) S_0}
  = \left\lVert \sum (1-U_0^*U_0) g_{l(n)} S_0 \right\rVert
  \le \sum \varepsilon\,2^{-n-1} = \varepsilon.
  \]
  By functional calculus, we can replace~\(U_0\) by another
  lifting~\(U\) of~\(\bar{U}\) such that \(US_0\) is exactly
  isometric and the range of \(U(1-S_0S_0^*)\) is orthogonal to
  the range of~\(US_0\).

  Next we construct~\(T\) and an isometry~\(S_1\) such that
  \(S\defeq S_0S_1\) and~\(T\) have all the required
  properties.  Let \((g_n)\) be as above and let \(h_n\defeq
  (US_0)g_n(US_0)^*\).  We have already observed that \((g_n)\)
  converges strictly to zero in \(\Mult(J)\).  So does~\(h_n\)
  because multiplication in \(\Mult(J)\) is separately strictly
  continuous.  Therefore, \(\norm{g_mh_n}\) converges to~\(0\)
  if either variable is fixed and the other goes to~\(\infty\).
  Thus we can find subsequences \(k(n)\) and \(l(n)\) such that
  \(\norm{g_{k(m)}h_{l(n)}}\le \varepsilon\,2^{-m-n-1}\).
  Lemma~\ref{lem:small_projections} provides isometries \(S_1\)
  and~\(T\) with
  \[
  \sum g_{l(n)}\cdot S_1 = S_1,
  \qquad
  \sum g_{k(n)}\cdot T = T.
  \]
  As above, the range of~\(T\) is orthogonal to~\(P\).
  Moreover,
  \[
  T^*US_0S_1 = T^* \sum_m g_{k(m)} US_0 \sum_n g_{l(n)} S_1 =
  T^* \sum_{m,n} g_{k(m)} h_{l(n)} US_0 S_1
  \]
  has norm at most \(\sum \varepsilon\,2^{-n-m-1} =
  \varepsilon\).  Thus the isometries \(S\defeq S_0S_1\)
  and~\(T\) have the desired properties.
\end{proof}

\section{Representable \texorpdfstring{$\K$}{K}-theory via crossed
  products}
\label{sec:alternative_RK}

As before, \(\Grd\) is a locally compact, second countable,
Hausdorff groupoid with Haar system and \(X\) and~\(Y\) are
proper, locally compact, second countable \(\Grd\)\nb-spaces.

Our next goal is to describe \(\RK^*_\Grd(X)\) and
\(\RK^*_{\Grd,Y}(X)\) as \(\K\)\nb-theory groups of certain
crossed products.  The relevant crossed products are not
\Cstar{}algebras but \sigCstar{}algebras
(see~\cite{Phillips:Inverse}), and we use the \(\K\)\nb-theory
for such algebras defined by Chris Phillips
in~\cite{Phillips:Representable_K}; he calls this theory
``representable \(\K\)\nb-theory'' to distinguish it from other
theories like algebraic \(\K\)\nb-theory, which may be
different even for~\(\K_0\).  We drop the adjective because
Phillips' theory is the only one with good homological
properties for such topological algebras.  The familiar
properties of \(\K\)\nb-theory for \Cstar{}algebras like Bott
periodicity, stability, homotopy invariance, and six-term exact
sequences for extensions all extend (see
\cite{Phillips:Representable_K}*{Theorem 3.4}).  We will use
this in \S\ref{sec:functorial} to study the homological
properties of the functor \(\RK^*_\Grd\).

\begin{definition}[see~\cite{Phillips:Inverse}]
  \label{def:sigma_Cstar}
  A \emph{\sigCstar{}algebra} is a complete topological
  \Star{}algebra whose topology is defined by an increasing
  sequence of \Cstar{}seminorms.  Equivalently, it is a
  countable projective limit of \Cstar{}algebras.
\end{definition}

The next theorem is a generalisation of
\cite{Emerson-Meyer:Descent}*{Lemma 20} and
\cite{Echterhoff-Emerson-Kim:Duality}*{Theorem 2.8} to
groupoids; the first result of this kind seems
\cite{Kasparov-Skandalis:Bolic}*{Theorem 5.4}.  For a better
perspective, we have stated a more general theorem than we
actually need, which involves \emph{bivariant} \(\K\)\nb-theory
for \sigCstar{}algebras.  Although it took some time to notice
this, Kasparov theory extends \emph{literally} to separable
\sigCstar{}algebras.  This was first carried out by Alexander
Bonkat in his thesis~\cite{Bonkat:Thesis}, using the extension
of Kasparov's Technical Lemma in~\cite{Hennings:Kasparov}.
Bonkat's definition is equivalent to the obvious one suggested
already by Chris Phillips in
\cite{Phillips:Representable_K}*{p.\ 470}.  The Kasparov group
\(\KK_*(\C,B)\) most relevant to us is defined in
\cite{Phillips:Representable_K}*{Definition 4.1} and identified
with Phillips' \(\K\)\nb-theory in
\cite{Phillips:Representable_K}*{Theorem 4.2}.

\begin{theorem}
  \label{the:RKK_via_crossed}
  Let \((X_n)\) be an increasing sequence of
  \(\Grd\)\nb-compact, \(\Grd\)\nb-invariant subsets of~\(X\)
  with \(\bigcup X_n=X\).  Let~\(A\) be a \Cstar{}algebra with
  the trivial action of~\(\Grd\) and let \(B\) be a
  \(\Grd\ltimes\nobreak X\)-\Cstar{}algebra.  Then there is a natural
  isomorphism
  \[
  \KK^{\Grd\ltimes X}_*(\CONT_0(X)\otimes A, B) \cong
  \KK_*(A,\varprojlim_n \Grd\ltimes B|_{X_n}).
  \]
  If~\(X\) is \(\Grd\)\nb-compact, then this agrees with \(\KK_*(A,
  \Grd\ltimes B)\).
\end{theorem}

Any \(\Grd\)\nb-compact, \(\Grd\)\nb-invariant subset of~\(X\)
is contained in~\(X_n\) for some \(n\in\N\), and the
subsets~\(X_n\) are automatically closed because the action is
proper.  Therefore, the restriction \(B|_{X_n}\) is simply the
quotient of~\(B\) by the ideal \(\CONT_0(X\setminus X_n)\cdot
B\).  The quotient maps turn \((\Grd\ltimes
B|_{X_n})_{n\in\N}\) into a countable projective system of
\Cstar{}algebras with surjective maps, whose projective limit
is a \sigCstar{}algebra.

All choices of the sequence \((X_n)\) are equivalent in a
suitable sense and therefore yield isomorphic
\sigCstar{}algebras \(\varprojlim_n \Grd\ltimes B|_{X_n}\).  In
particular, if~\(X\) itself is \(\Grd\)\nb-compact, then the
right hand side agrees with the usual Kasparov group
\(\KK_*(A,\Grd\ltimes B)\).

We mainly use the case \(B=\CONT_0(X)\), where we get
\(B|_{X_n} \cong \CONT_0(X_n)\).  The \sigCstar{}algebra
\(\varprojlim \CONT_0(X_n)\) consists of continuous functions
\(X\to\C\) whose restrictions to~\(X_n\) vanish at infinity for
all \(n\in\N\).

\begin{notation}
  \label{not:sigma_Cstar_crossed}
  We abbreviate
  \[
  \sC(\Grd\ltimes X) \defeq
  \varprojlim C^*(\Grd\ltimes X_n) \cong
  \varprojlim \Grd\ltimes \CONT_0(X_n).
  \]
\end{notation}

\begin{proof}[of Theorem~\ref{the:RKK_via_crossed}]
  Let \(\Grd\backslash X\) be the orbit space.
  \cite{Tu:Novikov}*{Proposition 6.25} yields
  \[
  \KK^{\Grd\ltimes X}_*(\CONT_0(X)\otimes A, B) \cong
  \KK^{\Grd\backslash X}_*(\CONT_0(\Grd\backslash X) \otimes A,
  \Grd\ltimes B).
  \]
  Combining this with the non-equivariant case of
  \cite{Emerson-Meyer:Descent}*{Lemma 20} yields the assertion.
  Alternatively, the proof of
  \cite{Emerson-Meyer:Descent}*{Lemma 20} extends to the
  groupoid case.
\end{proof}

\begin{corollary}
  \label{cor:RK_via_crossed}
  There is a natural isomorphism
  \[
  \RK^*_\Grd(X)
  \cong \K_*\bigl(\sC(\Grd\ltimes X)\bigr).
  \]
  If~\(X\) is \(\Grd\)\nb-compact, then \(\RK^*_\Grd(X) \cong
  \K_*\bigl(C^*(\Grd\ltimes X)\bigr) = \K^*_\Grd(X)\).
\end{corollary}

\begin{corollary}
  \label{cor:RK_support_via_crossed}
  Let \(p\colon X\to Y\) be a continuous
  \(\Grd\)\nb-equivariant map between two proper
  \(\Grd\)\nb-spaces, and let \((Y_n)_{n\in\N}\) be an
  increasing sequence of \(\Grd\)\nb-compact,
  \(\Grd\)\nb-invariant subspaces of~\(Y\).  Let \(X_n\defeq
  p^{-1}(Y_n)\).  Then there is a natural isomorphism
  \[
  \RK^*_{\Grd,Y}(X)
  \cong \K_*\bigl(\varprojlim \Grd\ltimes\CONT_0(X_n)\bigr)
  \cong \K_*\bigl(\varprojlim C^*(\Grd\ltimes X_n)\bigr).
  \]
\end{corollary}

The \sigCstar{}algebra \(\varprojlim C^*(\Grd\ltimes X_n)\) in
Corollary~\ref{cor:RK_support_via_crossed} may differ from
\(\sC(\Grd\ltimes X)\) because the subsets~\(X_n\) in
Corollary~\ref{cor:RK_support_via_crossed} need not be
\(\Grd\)\nb-compact.

Both corollaries follow directly from the definitions using
Theorem~\ref{the:RKK_via_crossed}.  The \Cstar{}algebras
\(\Grd\ltimes\CONT_0(X_n)\) and \(C^*(\Grd\ltimes X_n)\) are
isomorphic by~\eqref{eq:crossed_product_groupoid}.  We can also
derive the results in \S\ref{sec:Fred_compare_support} from
these corollaries.  For instance, we get once again that
\(\RK^*_{\Grd,Y}(X)=\K^*_\Grd(X)\) if~\(Y\) is
\(\Grd\)\nb-compact.

Phillips' representable \(\K\)\nb-theory for
\sigCstar{}algebras is computable by a Milnor
\(\varprojlim\nolimits^1\)-sequence
(\cite{Phillips:Representable_K}*{Theorem 3.2}).  This
specialises to short exact sequences
\begin{equation}
  \label{eq:Milnor_for_RK}
  \varprojlim\nolimits^1 \K^{*+1}_\Grd(X_n)
  \into \RK^*_{\Grd,Y}(X) \prto \varprojlim \K^*_\Grd(X_n),
\end{equation}
where \((X_n)\) is as in
Corollary~\ref{cor:RK_support_via_crossed}.  Thus we can, in
principle, reduce our theories to \(\K^*_\Grd(X_n)\).

\subsection{Functorial properties}
\label{sec:functorial}

First we discuss multiplication.  The Kasparov product turns
\[
\RK^*_\Grd(X) \defeq
\KK^{\Grd\ltimes X}\bigl(\CONT_0(X),\CONT_0(X)\bigr)
\]
into a graded-commutative ring with unit.  This ring acts by
exterior product on all Kasparov groups of the form
\(\KK^{\Grd\ltimes Z}(A,B)\) for \(\Grd\ltimes
X\)-\Cstar{}algebras \(A\) and~\(B\), where~\(Z\) is another
proper \(\Grd\)\nb-space with a map \(p\colon X\to Z\); we view
\(\Grd\ltimes X\)-\Cstar{}algebras as \(\Grd\ltimes
Z\)-\Cstar{}algebras via the resulting forgetful map.  In
particular, the ring \(\RK^*_\Grd(X)\) acts on
\(\RK^*_{\Grd,Y}(X)\), that is, the latter is a graded module
over \(\RK^*_\Grd(X)\) in a canonical way.  Furthermore, via
the descent homomorphism
\begin{multline*}
  \RK^*_\Grd(X) \defeq \KK^{\Grd\ltimes
    X}\bigl(\CONT_0(X),\CONT_0(X)\bigr)
  \\ \xrightarrow{\textup{descent}}
  \KK\bigl(\Grd\ltimes \CONT_0(X), \Grd\ltimes\CONT_0(X)\bigr)
  = \KK\bigl(C^*(\Grd\ltimes X), C^*(\Grd\ltimes X)\bigr),
\end{multline*}
the ring \(\RK^*_\Grd(X)\) acts on \(\K^*_\Grd(X) \defeq
\K_*\bigl(C^*(\Grd\ltimes X)\bigr)\).  Thus \(\K^*_\Grd(X)\) is
a graded module over the graded ring \(\RK^*_\Grd(X)\).  Via
the canonical map \(\K^*_\Grd(X)\to \RK^*_\Grd(X)\), we may
also equip \(\K^*_\Grd(X)\) with a ring structure in its own
right, but this ring need not be unital.

\smallskip

Now we turn to functoriality.  We only discuss
\(\RK^*_\Grd(X)\) and \(\K^*_\Grd(X)\).

A continuous \(\Grd\)\nb-equivariant map \(f\colon X\to Y\)
induces a grading preserving ring homomorphism
\[
\RK^*_\Grd(f)\colon \RK^*_\Grd(Y)\to \RK^*_\Grd(X),
\]
regardless whether~\(f\) is proper.  Hence \(X\mapsto
\RK^*_\Grd(X)\) is a contravariant functor from the category of
locally compact proper \(\Grd\)\nb-spaces to the category of
graded-commutative graded rings.  Similarly, a continuous,
\emph{proper}, \(\Grd\)\nb-equivariant map induces a grading
preserving group homomorphism
\[
\K^*_\Grd(f)\colon \K^*_\Grd(Y)\to \K^*_\Grd(X).
\]
Both functors \(\RK^*_\Grd\) and \(\K^*_\Grd\) are homotopy
invariant (the latter only for proper homotopies, of course)
because \(\K\)\nb-theory for \sigCstar{}algebras is homotopy
invariant.

More generally, any continuous groupoid homomorphism (functor)
\(f\colon \Grd\ltimes X\to \Grd'\ltimes X'\) induces maps
\[
f^*\colon \KK^{\Grd'\ltimes X'}(A,B) \to \KK^{\Grd\ltimes X}(f^*A,f^*B)
\]
by \cite{LeGall:KK_groupoid}*{\S7.1}, where \(f^*A \defeq
\CONT_0(X)\otimes_{\CONT_0(X')} A\) equipped with the canonical
action of~\(\Grd\).  These are compatible with Kasparov
products and satisfy the functoriality properties \((f\circ
g)^*=g^*\circ f^*\) and \(\ID^*=\ID\) as expected.  Since
\(f^*\CONT_0(X') \cong \CONT_0(X)\), this specialises to a ring
homomorphism \(\RK^*_{\Grd'}(X')\to \RK^*_\Grd(X)\), that is,
equivariant representable \(\K\)\nb-theory is functorial for
strict groupoid morphisms.

\begin{theorem}
  \label{the:K_RK_Morita_invariant}
  A Morita equivalence between the two groupoids \(\Grd\ltimes
  X\) and \(\Grd'\ltimes X'\) induces isomorphisms
  \(\RK^*_\Grd(X) \cong \RK^*_{\Grd'}(X')\) and \(\K^*_\Grd(X)
  \cong \K^*_{\Grd'}(X')\).
\end{theorem}

\begin{proof}
  Morita equivalent groupoids give rise to equivalent Kasparov
  categories by \cite{LeGall:KK_groupoid}*{\S7.2}.  In
  particular, this yields the isomorphism \(\RK^*_\Grd(X) \cong
  \RK^*_{\Grd'}(X')\).  Furthermore, Morita equivalent
  groupoids have Morita equivalent \Cstar{}algebras
  by~\cite{Muhly-Renault-Williams:Equivalence}.  Therefore,
  \(\K^*_\Grd(X) \cong \K^*_{\Grd'}(X')\) as well.
\end{proof}

\begin{example}
  \label{exa:induction}
  Let~\(G\) be a group and let~\(H\) be a closed subgroup
  in~\(G\).  Then \(G\ltimes G/H\) is Morita equivalent
  to~\(H\).  Therefore, we have canonical induction
  isomorphisms
  \[
  \RK^*_G(G\times_H X) \cong \RK^*_H(X),
  \qquad
  \K^*_G(G\times_H X) \cong \K^*_H(X)
  \]
  for any \(H\)\nb-space~\(X\).
\end{example}

\begin{example}
  \label{exa:free_action}
  Let~\(X\) be a free and proper \(\Grd\)\nb-space.  Then
  \(\Grd\ltimes X\) is Morita equivalent to the orbit space
  \(\Grd\backslash X\).  Hence
  \[
  \RK^*_\Grd(X) \cong \RK^*(\Grd\backslash X),
  \qquad
  \K^*_\Grd(X) \cong \K^*(\Grd\backslash X).
  \]
\end{example}

\begin{proposition}
  \label{pro:K_RK_additive}
  Let \((X_n)_{n\in\N}\) be \(\Grd\)\nb-spaces and let
  \(X\defeq \bigsqcup_{n\in\N} X_n\) be their disjoint union.
  Then
  \[
  \K^*_\Grd(X) \cong \bigoplus_{n\in\N} \K^*_\Grd(X_n),
  \qquad
  \RK^*_\Grd(X) \cong \prod_{n\in\N} \RK^*_\Grd(X_n).
  \]
\end{proposition}

\begin{proof}
  The first isomorphism follows from the additivity of
  \(\K\)\nb-theory for direct sums of \Cstar{}algebras.  The
  second isomorphism follows from the behaviour of
  \(\K\)\nb-theory for \sigCstar{}algebras for direct products
  (see \cite{Phillips:Representable_K}*{Proposition 3.1}).
\end{proof}

\smallskip

Finally, we turn to exact sequences.  To formulate them, we
extend our theories to pairs of spaces \((X,A)\), where~\(X\)
is a proper \(\Grd\)\nb-space and \(A\subseteq X\) is a closed
\(\Grd\)\nb-invariant subspace.

For \(\K^*_\Grd\), we simply put
\begin{equation}
  \label{eq:KG_definition}
  \K^*_\Grd(X,A) \defeq \K^*_\Grd(X\setminus A)
  = \K_*\bigl(\Grd \ltimes \CONT_0(X\setminus A)\bigr).
\end{equation}
This extends the old theory because \(\K^*_\Grd(X) =
\K^*_\Grd(X,\emptyset)\).  The excision property
\begin{equation}
  \label{eq:KG_excision}
  \K^*_\Grd(X,A) \cong \K^*_\Grd(X\setminus U,A\setminus U)
\end{equation}
for an open, \(\Grd\)\nb-invariant subset \(U\subseteq A\) is
trivial.  Since the (full) crossed product functor is exact, we
get an extension of \Cstar{}algebras
\[
\Grd\ltimes \CONT_0(X\setminus A) \into
\Grd\ltimes \CONT_0(X) \prto \Grd\ltimes \CONT_0(A)
\]
for any pair \((X,A)\), leading to an exact sequence
\begin{equation}
  \label{eq:KG_exact}
  \begin{gathered}
    \xymatrix{
      \K_\Grd^0(X,A) \ar[r] &
      \K_\Grd^0(X) \ar[r] &
      \K_\Grd^0(A) \ar[d] \\
      \K_\Grd^1(A) \ar[u] &
      \K_\Grd^1(X) \ar[l] &
      \K_\Grd^1(X,A). \ar[l]
    }
  \end{gathered}
\end{equation}
Pairs of \(\Grd\)\nb-spaces form a category, whose morphisms
are \(\Grd\)\nb-equivariant continuous proper maps that
restrict to maps between the specified subspaces.  It is clear
that \(\K^*_\Grd\) is a contravariant functor on this category.
The isomorphism in~\eqref{eq:KG_excision} is induced by the
obvious morphism \((X\setminus U,A\setminus U)\to (X,A)\), the
horizontal maps in~\eqref{eq:KG_exact} are induced by the
obvious morphisms \((A,\emptyset)\to (X,\emptyset) \to (X,A)\).

For \(\RK^*_\Grd\), we put
\begin{equation}
  \label{eq:RK_G_definition}
  \RK^*_\Grd(X,A) \defeq \RK^*_{\Grd,X}(X\setminus A).
\end{equation}
Corollary~\ref{cor:RK_support_via_crossed} identifies this with
\(\K_*\bigl(\varprojlim \Grd\ltimes \CONT_0(X_n\setminus
A)\bigr)\), where \((X_n)\) is an increasing sequence of
\(\Grd\)\nb-compact \(\Grd\)\nb-invariant subsets of~\(X\) with
\(\bigcup X_n = X\).  Again, we have \(\RK^*_\Grd(X) =
\RK^*_\Grd(X,\emptyset)\).  The excision property
\begin{equation}
  \label{eq:RKG_excision}
  \RK^*_\Grd(X,A) \cong \RK^*_\Grd(X\setminus U,A\setminus U)
\end{equation}
holds for any open, \(\Grd\)\nb-invariant subset \(U\subseteq
A\) because both sides are computed by the same
\sigCstar{}algebra.  As above, we get a sequence of
\Cstar{}algebra extensions
\[
\Grd\ltimes \CONT_0(X_n\setminus A) \into
\Grd\ltimes \CONT_0(X_n) \prto
\Grd\ltimes \CONT_0(X_n\cap A)
\]
for all \(n\in\N\).  The vertical maps
\[
\xymatrix{
  \Grd\ltimes \CONT_0(X_{n+1}\setminus A) \ar[r] \ar[d] &
  \Grd\ltimes \CONT_0(X_{n+1}) \ar[r] \ar[d] &
  \Grd\ltimes \CONT_0(X_{n+1}\cap A) \ar[d] \\
  \Grd\ltimes \CONT_0(X_n\setminus A) \ar[r] &
  \Grd\ltimes \CONT_0(X_n) \ar[r] &
  \Grd\ltimes \CONT_0(X_n\cap A)
}
\]
are all surjective.  Therefore, we get an extension of
\sigCstar{}algebras
\[
\varprojlim \Grd\ltimes \CONT_0(X_n\setminus A) \into
\sC(\Grd\ltimes X_n) \prto
\sC(\Grd\ltimes X_n\cap A),
\]
which induces a six-term exact sequence
\begin{equation}
  \label{eq:RKG_exact}
  \begin{gathered}
    \xymatrix{
      \RK_\Grd^0(X,A) \ar[r] &
      \RK_\Grd^0(X) \ar[r] &
      \RK_\Grd^0(A) \ar[d] \\
      \RK_\Grd^1(A) \ar[u] &
      \RK_\Grd^1(X) \ar[l] &
      \RK_\Grd^1(X,A) \ar[l]
    }
  \end{gathered}
\end{equation}
by \cite{Phillips:Representable_K}*{Theorem 3.4}.  A
\(\Grd\)\nb-equivariant continuous map of pairs \((X,A)\to
(X',A')\) induces a map \(\RK^*_\Grd(X',A') \to
\RK^*_\Grd(X,A)\) even if it is not proper.  The horizontal
maps in~\eqref{eq:RKG_exact} and the excision isomorphism
in~\eqref{eq:RKG_excision} can be described as for
\(\K^*_\Grd\).

We define a relative version of \(\RK^*_{\Grd,Y}\) by
\[
\RK^*_{\Grd,Y}(X,A) \defeq \RK^*_{\Grd,Y}(X\setminus A)
\]
for a pair \((X,A)\) of \(\Grd\ltimes Y\)\nb-spaces.  Like
\(\K^*_\Grd\), this theory is functorial for proper morphisms
of pairs, satisfies excision, and has long exact sequences.

The natural transformations \(\K^*_\Grd(X)\to \RK^*_{\Grd,Y}(X)\to
\RK^*_\Grd(X)\) extend to natural transformations
\[
\K^*_\Grd(X,A) \to \RK^*_{\Grd,Y}(X,A) \to \RK^*_\Grd(X,A).
\]
The first of these maps is invertible if~\(Y\) is
\(\Grd\)\nb-compact, the second is invertible if the closure of
\(X\setminus A\) in~\(X\) is \(Y\)\nb-compact: in both cases,
the assumption ensures that both sides are computed by the same
\sigCstar{}algebra.

\begin{lemma}
  \label{lem:K_support_relative}
  Any \(\Grd\)\nb-space~\(X\) over~\(Y\) is of the form
  \(X=\bar{X}\setminus\partial X\) for some \(Y\)\nb-compact
  \(\Grd\)\nb-space~\(\bar{X}\) and some closed
  \(\Grd\)\nb-invariant subset~\(\partial X\) of~\(\bar{X}\),
  necessarily \(Y\)\nb-compact.  We have
  \[
  \RK^*_{\Grd,Y}(X) \cong \RK_\Grd^*(\bar{X},\partial X).
  \]
\end{lemma}

\begin{proof}
  The second assertion is just~\eqref{eq:RK_G_definition}, so
  that it remains to construct \(\bar{X}\) and~\(\partial X\).
  Roughly speaking, we construct~\(\bar{X}\) by one-point
  compactifying the fibres \({p^{-1}(y)\subseteq X}\) of the
  map \(p\colon X\to Y\) for all \({y\in Y}\).  More precisely,
  we let \(\bar{X}\defeq X\sqcup Y\) with the unique topology
  where~\(X\) is open and carries the given topology, and
  \(U_X\sqcup U_Y\subseteq X\sqcup Y\) is a neighbourhood of
  \(y\in Y\) if and only if~\(U_Y\) is a neighbourhood of~\(y\)
  in~\(Y\) and there is a neighbourhood~\(V\) of~\(y\) in~\(Y\)
  such that \(p^{-1}(V)\setminus U_X\) is relatively compact.
  It is left to the reader to check that~\(\bar{X}\) is a
  locally compact Hausdorff space, that the map \(\bar{X}\to
  Y\) is continuous and proper, and that the \(\Grd\)\nb-action
  on~\(\bar{X}\) is continuous and proper.
\end{proof}

For the special choice of \(\bar{X}\) and \(\partial X\) in the
proof of Lemma~\ref{lem:K_support_relative}, \(\partial X=Y\)
is a retract of~\(\bar{X}\) via the projection map \(\bar{X}\to
Y\).  Hence the extension of \sigCstar{}algebras
\[
\varprojlim \Grd\ltimes \CONT_0(X_n) \into
\sC(\Grd\ltimes \bar{X}_n) \prto
\sC(\Grd\ltimes \partial X_n)
\]
splits.  The resulting \(\K\)\nb-theory long exact sequence
splits as well, that is,
\begin{equation}
  \label{eq:K_support_as_kernel}
  \RK_{\Grd,Y}^*(X) \cong \RK_\Grd^*(\bar{X},\partial X)
  \cong \ker\bigl(\RK_\Grd^*(\bar{X})
  \to \RK_\Grd^*(\partial X)\bigr).
\end{equation}
This is analogous to the familiar description of~\(\K_0\) for
non-unital \Cstar{}algebras as \(\K_0(A)\defeq \ker
\bigl(\K_0(A^+) \to \K_0(\C)\bigr)\).

\begin{theorem}
  \label{the:RK_support_continuous}
  Let~\(X\) be a space over~\(Y\) and let~\(I_Y\) be the
  directed set of all relatively \(Y\)\nb-compact,
  \(\Grd\)\nb-invariant, open subsets of~\(X\) \textup{(}that
  is, sets in~\(I_Y\) have \(Y\)\nb-compact closure
  in~\(X\)\textup{)}.  Then
  \[
  \RK^*_{\Grd,Y}(X)
  \cong \varinjlim_{A\in I_Y} \RK^*_{\Grd,Y}(A)
  \cong \varinjlim_{A\in I_Y} \RK^*_\Grd(X,X\setminus A).
  \]
  Similarly, let~\(I_\Grd\) be the directed set of all
  relatively \(\Grd\)\nb-compact, \(\Grd\)\nb-invariant, open
  subsets of~\(X\).  Then
  \[
  \K^*_\Grd(X)
  \cong \varinjlim_{A\in I_\Grd} \K^*_\Grd(A)
  \cong \varinjlim_{A\in I_\Grd} \RK^*_\Grd(X,X\setminus A).
  \]
\end{theorem}

\begin{proof}
  The isomorphism \(\K^*_\Grd(X) \cong \varinjlim_{A\in I_\Grd}
  \K^*_\Grd(A)\) follows because
  \[
  \Grd\ltimes \CONT_0(X) \cong
  \varinjlim_{A\in I_\Grd}\Grd\ltimes \CONT_0(A)
  \]
  and \(\K\)\nb-theory for \Cstar{}algebras is continuous with
  respect to direct limits.  If \(A\in I_\Grd\),
  then~\eqref{eq:KG_definition}, the \(\Grd\)\nb-compactness
  of~\(\cl{A}\), and~\eqref{eq:RKG_excision} yield
  \[
  \K^*_\Grd(A) \cong \K^*_\Grd(\cl{A},\cl{A}\setminus A)
  \cong \RK^*_\Grd(\cl{A},\cl{A}\setminus A)
  \cong \RK^*_\Grd(X,X\setminus A).
  \]

  We have an isomorphism
  \[
  \varprojlim_{n\in\N} \Grd\ltimes \CONT_0(X_n) \cong
  \varinjlim_{A\in I_Y} \varprojlim_{n\in\N}
  \Grd\ltimes \CONT_0(A\cap X_n)
  \]
  as well, where the subsets \(X_n\subseteq X\) are defined as
  in Corollary~\ref{cor:RK_support_via_crossed}.  But inductive
  limits of \sigCstar{}algebras require some care because they
  may become uncountable projective systems of
  \Cstar{}algebras.  Hence we prefer another argument in this
  case.

  Theorem~\ref{the:RK_Fred_general} identifies
  \(\RK^*_{\Grd,Y}(X)\) with the set of homotopy classes of
  continuous \(\Grd\)\nb-equivariant maps \(U\colon
  X\to\Fred_\Grd\) with \(Y\)\nb-compact support.  Any
  \(Y\)\nb-compact subset of~\(X\) has a relatively
  \(Y\)\nb-compact open neighbourhood.  Let~\(A\) be such a
  neighbourhood of the support of~\(U\).  The restriction
  of~\(U\) to~\(A\) is a cycle for \(\RK^*_{\Grd,Y}(A)\).
  Similarly, if two cycles for \(\RK^*_{\Grd,Y}(X)\) are
  homotopic, then we can find a relatively \(Y\)\nb-compact
  neighbourhood~\(B\) for the support of the resulting homotopy
  and get a homotopy in \(\RK^*_{\Grd,Y}(B)\).  It remains to
  check that the canonical map
  \(\RK^*_{\Grd,Y}(A)\to\RK^*_{\Grd,Y}(X)\) maps these
  restricted cycles to the original cycles.  We must describe,
  therefore, how the latter map acts on maps
  \(A\to\Fred_\Grd\).

  Let \(U\colon A\to\Fred_\Grd\) be a \(\Grd\)\nb-equivariant
  continuous map whose support is closed in~\(X\).  View~\(U\)
  as a \(\Grd\)\nb-equivariant adjointable operator on the
  Hilbert module \(L^2(\Grd\ltimes A)^\infty\) over
  \(\CONT_0(A)\).  Since \(\CONT_0(A)\subseteq \CONT_0(X)\) is
  an ideal, we may view \(L^2(\Grd\ltimes A)^\infty\) as a
  Hilbert module over \(\CONT_0(X)\) and~\(U\) as an
  adjointable operator of Hilbert \(\CONT_0(X)\)\nb-modules.
  If~\(U\) is the restriction of a map \(\cl{U}\colon
  X\to\nobreak\Fred_\Grd\) with the same support, then the
  cycle defined by~\(\cl{U}\) on \(L^2(\Grd\ltimes X)^\infty\)
  is homotopic to the one defined by~\(U\) on \(L^2(\Grd\ltimes
  A)^\infty\): the homotopy is given by the
  \(\Grd\)\nb-equivariant Hilbert \(\CONT_0(X)\)-module
  \[
  \{f\in \CONT([0,1],L^2(\Grd\ltimes X)^\infty) \mid
  f(0)\in L^2(\Grd\ltimes A)^\infty\}
  \]
  with the restriction of~\(\cl{U}\).  This is indeed a cycle
  for \(\RK^0_{\Grd,Y}([0,1]\times X)\) because~\(\cl{U}\) is
  invertible on the closure of \(X\setminus A\).

  The proof that \(\RK^*_{\Grd,Y}(X) \cong \varinjlim_{A\in
    I_Y} \RK^*_{\Grd,Y}(A)\) is now finished easily.
  Since~\(\cl{A}\) is \(Y\)\nb-compact for all \(A\in I_Y\),
  Lemma~\ref{lem:K_support_relative} and the excision property
  of \(\RK^*_\Grd\) yield
  \[
  \RK^*_{\Grd,Y}(A) \cong
  \RK^*_\Grd(\cl{A},\cl{A}\setminus A) \cong
  \RK^*_\Grd(X,X\setminus A).\qedhere
  \]
\end{proof}

\section{Equivariant vector bundles and the crossed product}
\label{sec:RK_via_vector_bundle}

We are going to relate \(\Grd\)\nb-equivariant vector bundles
over~\(X\) to the \sigCstar{}algebra \(\sC(\Grd\ltimes X)\).
As before, \(\Grd\) is a second countable, locally compact,
Hausdorff groupoid with Haar system and~\(X\) is a second
countable, locally compact, proper \(\Grd\)\nb-space.

We fix an increasing sequence \((X_n)_{n\in\N}\) of
\(\Grd\)\nb-compact \(\Grd\)\nb-invariant subsets with
\(\bigcup X_n=X\) and abbreviate
\[
B_n\defeq C^*(\Grd\ltimes X_n) \cong \Grd\ltimes \CONT_0(X_n),
\qquad
B\defeq \varprojlim B_n
= \sC(\Grd\ltimes X_n).
\]
Recall that a Hilbert module over~\(B\) is of the form \(\Hilm
= \varprojlim \Hilm_n\), where each~\(\Hilm_n\) is a Hilbert
module over~\(B_n\) and we form the limit with respect to
\(B_{n+1}\)\nb-linear projections \(\Hilm_{n+1}\to\Hilm_n\)
that induce unitary operators \(\Hilm_{n+1}\otimes_{B_{n+1}}
B_n\cong \Hilm_n\).

The \sigCstar{}algebras of adjointable and compact operators
on~\(\Hilm\) are
\[
\Bound(\Hilm) = \varprojlim \Bound(\Hilm_n)
\quad\text{and}\quad
\Comp(\Hilm) = \varprojlim \Comp(\Hilm_n),
\]
respectively, where the inverse limits are taken in the
category of \sigCstar{}algebras and hence contain unbounded
elements.  Let \(\Bound\bound(\Hilm)\subseteq \Bound(\Hilm)\)
and \(\Comp\bound(\Hilm)\subseteq \Comp(\Hilm)\) be the
\Cstar{}algebras of bounded elements, that is, elements
\((T_n)_{n\in\N}\) with \(\sup {}\norm{T_n} <\infty\).

Let~\(\Cat_B\) be the category of Hilbert modules over~\(B\)
with \emph{bounded} adjointable operators as morphisms.  Let
\(\Cat_{\Grd,X}\) be the category of \(\Grd\)\nb-equivariant
Hilbert modules over \(\CONT_0(X)\), with
\(\Grd\)\nb-equivariant adjointable operators as morphisms.
Both \(\Cat_B\) and~\(\Cat_{\Grd,X}\) are evidently
\Cstar{}categories.

\begin{notation}
  \label{def:stabilisation}
  Let~\(B_\Comp\) for a \sigCstar{}algebra \(B=\varprojlim
  B_n\) be its stabilisation
  \[
  B_\Comp \defeq B\otimes\Comp(\ell^2\N)
  \defeq \varprojlim B_n\otimes \Comp(\ell^2\N).
  \]
\end{notation}

\begin{theorem}
  \label{the:category_equivalence_equivariant_crossed}
  The \Cstar{}categories \(\Cat_B\) and~\(\Cat_{\Grd,X}\) are
  equivalent.

  A \(\Grd\)\nb-equivariant Hilbert \(\CONT_0(X)\)\nb-module is
  the space of \(\CONT_0\)\nb-sections of a
  \(\Grd\)\nb-equivariant Hermitian \emph{vector bundle}
  over~\(X\) if and only if the associated Hilbert
  \(B\)\nb-module~\(\Hilm\) satisfies the following equivalent
  conditions:
  \begin{enumerate}[label=\textup{(\arabic{*})}]
  \item \(\Comp\bound(\Hilm)= \Bound\bound(\Hilm)\);
  \item \(\Comp(\Hilm)= \Bound(\Hilm)\);
  \item \(\ID_\Hilm\in\Comp\bound(\Hilm)\);
  \item \(\ID_\Hilm\in\Comp(\Hilm)\);
  \item there is a projection \(p\in B_\Comp\) whose range is
    isomorphic to~\(\Hilm\).
  \end{enumerate}
\end{theorem}

\begin{proof}
  The first assertion is implicit in the proof of
  \cite{Emerson-Meyer:Descent}*{Lemma 20} (in the group case).
  A \(\Grd\)\nb-equivariant Hilbert
  \(\CONT_0(X)\)-module~\(\HilmF\) restricts to
  \(\Grd\)\nb-equivariant Hilbert
  \(\CONT_0(X_n)\)-modules~\(\HilmF_n\) for all \(n\in\N\).
  Since~\(X_n\) is \(\Grd\)\nb-compact, these correspond to
  Hilbert \(B_n\)\nb-modules~\(\Hilm_n\) by the construction
  before \cite{Tu:Novikov}*{Proposition 6.24}, and this passage
  identifies \(\Grd\)\nb-equivariant adjointable operators
  on~\(\HilmF_n\) and adjointable operators on~\(\Hilm_n\).
  Moreover, compact operators on~\(\Hilm_n\) correspond to
  \(\Grd\)\nb-equivariant adjointable operators \(T\colon
  \HilmF_n\to\HilmF_n\) with \(M_f T\in\Comp(\HilmF_n)\) for
  all \(f\in\CONT_0(X_n)\), by \cite{Tu:Novikov}*{Proposition
    6.24}.

  Now we form the Hilbert \(B\)\nb-module \(\Hilm\defeq
  \varprojlim \Hilm_n\).  An adjointable operator on~\(\Hilm\)
  restricts to adjointable operators on~\(\Hilm_n\) for all
  \(n\in\N\), and, conversely, a bounded, compatible sequence
  of adjointable operators on~\(\Hilm_n\) must come from an
  adjointable operator on~\(\Hilm\).  This reduces the
  isomorphism \(\Bound\bound(\Hilm)\cong \Bound(\HilmF)^\Grd\)
  to the already know isomorphisms \(\Bound(\Hilm_n) \cong
  \Bound(\HilmF_n)^\Grd\) and yields the asserted equivalence
  of categories.

  We can also describe the compact operators:
  \begin{align*}
    \Comp\bound(\Hilm) &\cong \{T\in\Bound(\HilmF)^\Grd \mid
    \text{\(\varphi(f)T|_{X_n}\in\Comp(\HilmF_n)\) for all \(n\in\N\),
      \(f\in\CONT_0(X_n)\)}\}\\
    &\cong\{T\in\Bound(\HilmF)^\Grd \mid
    \text{\(M_fT\in\Comp(\HilmF)\) for all \(f\in\CONT_0(X)\)}\}.
  \end{align*}

  The Hilbert module~\(\HilmF\) is the space of
  \(\CONT_0\)\nb-sections of a vector bundle if and only if
  \(M_f\in\Comp(\HilmF)\) for all \(f\in\CONT_0(X)\); here we
  may disregard the \(\Grd\)\nb-action.  For the proof,
  write~\(\HilmF\) as a direct summand of \(\CONT_0(X)^\infty\)
  and let~\(p\) be the projection onto~\(\HilmF\).  Since \(M_f
  p\in\Comp(\CONT_0(X)^\infty) \cong \CONT_0(X,\Comp)\), this
  must be a norm-continuous map from~\(X\) to the space of
  projections in \(\Comp\defeq \Comp(\ell^2\N)\).  Norm
  continuity implies that the rank of the projection~\(p\) is
  locally constant, so that the range of~\(p\) is a vector
  bundle.

  As a result, \(\Hilm\) corresponds to the space of
  \(\CONT_0\)\nb-sections of a \(\Grd\)\nb-equivariant vector
  bundle \emph{if and only if}
  \(\ID_\Hilm\in\Comp\bound(\Hilm)\).  It remains to check that
  the five conditions on~\(\Hilm\) in the theorem are indeed
  equivalent.  Since \(\ID_\Hilm\) is bounded,
  \(\ID_\Hilm\in\Comp\bound(\Hilm)\) is equivalent to
  \(\ID_\Hilm\in\Comp(\Hilm)\).  Since \(\Comp\bound(\Hilm)\)
  and \(\Comp(\Hilm)\) are ideals in \(\Bound\bound(\Hilm)\)
  and \(\Bound(\Hilm)\), respectively, these are equivalent to
  \(\Comp\bound(\Hilm)= \Bound\bound(\Hilm)\) and
  \(\Comp(\Hilm)= \Bound(\Hilm)\).

  By the Stabilisation Theorem for Hilbert modules over
  \sigCstar{}algebras, which is proved like its \Cstar{}algebra
  counterpart (see~\cite{Mingo-Phillips:Triviality}), for any
  Hilbert module~\(\Hilm\) over~\(B\) there is an adjointable
  isometry \(S\colon \Hilm\to B^\infty\defeq
  B\otimes\ell^2(\N)\).  Thus~\(\Hilm\) is isomorphic to the
  range of the projection \(p\defeq SS^*\).  If \(\ID_\Hilm\)
  is compact, then so is \(S=S\circ\ID_\Hilm\) and hence
  \(p=S\circ S^*\).  Conversely, if~\(p\) is compact, so is
  \(S=p\circ S\) and hence \(\ID_\Hilm=S^*\circ S\).  Hence all
  five conditions on~\(\Hilm\) in the statement of the theorem
  are equivalent as asserted.
\end{proof}

\begin{corollary}
  \label{cor:equivariant_vector_bundle_as_K00}
  The monoid of isomorphism classes of \(\Grd\)\nb-equivariant
  \textup{(}Hermitian\textup{)} vector bundles on~\(X\) is
  isomorphic to the monoid of idempotents \textup{(}or
  projections\textup{)} in the stable \sigCstar{}algebra
  \(B_\Comp\defeq \sC(\Grd\ltimes X)_\Comp\).
\end{corollary}

\begin{proof}
  Theorem~\ref{the:category_equivalence_equivariant_crossed}
  shows that the monoid of isomorphism classes of Hermitian
  \(\Grd\)\nb-equivariant vector bundles on~\(X\) is equivalent
  to the monoid of equivalence classes of projections
  in~\(B_\Comp\).

  We claim that any \(\Grd\)\nb-equivariant vector bundle~\(V\)
  carries a \(\Grd\)\nb-equivariant Hermitian structure.  Its
  construction uses three ingredients.  First, we need a
  Hermitian structure \(({\cdot},{\cdot})\) on~\(V\) that need
  not be \(\Grd\)\nb-equivariant -- this exists because~\(X\)
  is second countable and locally compact, hence paracompact.
  Secondly, we need a Haar system on~\(\Grd\), that is, a left
  invariant continuous family \((\mu_\base)_{\base\in\Base}\)
  of non-negative measures on the fibres of the range map
  \(\Grd\to\Base\).  Thirdly, we need a cut-off function, that
  is, a function \(\varphi\colon X\to[0,\infty)\) with
  \(\Grd\)\nb-compact support such that \(\int \varphi(g\cdot
  x) \,\textup d\mu(g) = 1\) for all \(x\in X\) -- this exists
  because~\(\Grd\) acts properly and \(\Grd\backslash X\) is
  paracompact.  We define
  \[
  \langle\xi,\eta\rangle \defeq
  \int \varphi(g^{-1} x) (g^{-1}\xi,g^{-1}\eta)
  \,\textup d\mu(g)
  \qquad\text{for \(x\in X\), \(\xi,\eta\in V_x\).}
  \]
  This is the desired \(\Grd\)\nb-invariant Hermitian inner
  product.

  Furthermore, any two \(\Grd\)\nb-invariant Hermitian inner
  products on~\(V\) are homotopic by an affine homotopy.  Hence
  it makes no difference whether we study vector bundles or
  Hermitian vector bundles.  Moreover, any idempotent in a
  \sigCstar{}algebra is equivalent to a projection, that is, a
  self-adjoint idempotent, so that it makes no difference
  whether we use idempotents or projections.
\end{proof}

Since projections are automatically bounded, we may
replace~\(B_\Comp\) with its \Cstar{}subalgebra of bounded
elements in
Corollary~\ref{cor:equivariant_vector_bundle_as_K00}.  But
there are good reasons to work with \sigCstar{}algebras in this
context, such as the non-separability of \((B_\Comp)\bound\)
and the topology on unitary groups (see
\cite{Phillips:Representable_K}*{Example 1.2}).

Now we can express \(\vbK^0_\Grd(X)\) and, more generally, the
relative theory \(\vbK^0_\Grd(X,A)\), in terms of crossed
products:

\begin{definition}
  \label{def:Vect_B}
  Let~\(B\) be a \sigCstar{}algebra.  We let \(\Vect(B)\) be
  the set of Murray--von Neumann equivalence classes of
  projections in the stabilisation~\(B_\Comp\).  This is a
  monoid with respect to the orthogonal direct sum, where we
  use a canonical isomorphism \(\Mat_2(B_\Comp)\cong B_\Comp\).
  We let \(\K_{00}(B)\) the Grothendieck group of \(\Vect(B)\).
\end{definition}

\begin{definition}
  \label{def:Vect_B_rel}
  Let~\(B\) be a \sigCstar{}algebra, \(J\subseteq B\) a closed
  \Star{}ideal, and \(\pi\colon B_\Comp\to (B/J)_\Comp\) the
  quotient map.  Consider triples \((p^+,p^-,v)\) consisting of
  two projections \(p^\pm\in B_\Comp\) and a partial isometry
  \(v\in (B/J)_\Comp\) whose range and source projections are
  \(\pi(p^-)\) and~\(\pi(p^+)\), respectively, that is,
  \(vv^*=\pi(p^-)\), \(v^*v=\pi(p^+)\).  Two such triples
  \((p_t^+,p_t^-,v_t)\) for \(t=0,1\) are considered
  \emph{equivalent} if there are partial isometries \(s^+,s^-\)
  in~\(B_\Comp\) with
  \[
  s^\pm (s^\pm)^* = p_1^\pm,\qquad
  (s^\pm)^* s^\pm = p_0^\pm,\qquad
  \pi(s^-) v_0 = v_1 \pi(s^+).
  \]
  The direct sum of such triples is defined as usual and turns
  the set of equivalence classes of triples into a monoid,
  which we denote by \(\Vect(B,B/J)\).

  Call a triple \((p^+,p^-,v)\) \emph{degenerate} if there is a
  partial isometry \(\hat{v}\in B_\Comp\) with
  \(\pi(\hat{v})=v\), \(\hat{v}\hat{v}^*=p^-\), and
  \(\hat{v}^*\hat{v}=p^+\).  The degenerate triples form a
  submonoid in \(\Vect(B,B/J)\).  Let \(\K_{00}(B,B/J)\) be the
  quotient of \(\Vect(B,B/J)\) by addition of degenerate
  cycles.
\end{definition}

It is easy to see that \(\K_{00}(B,B/J)\) is always a group.
The unit element is the class of the degenerate triples, the
inverse of \((p^+,p^-,v)\) is \((p^-,p^+,v^*)\).

\begin{theorem}
  \label{the:vbK_via_crossed}
  Let~\(X\) be a proper \(\Grd\)\nb-space and let \(A\subseteq
  X\) be a closed \(\Grd\)\nb-invariant subset.  Define
  \[
  B\defeq \sC(\Grd\ltimes X),
  \qquad
  J \defeq \varprojlim \Grd\ltimes \CONT_0(X_n\setminus A).
  \]
  Then~\(J\) is an ideal in~\(B\).  There are canonical monoid
  isomorphisms
  \[
  \Vect_\Grd(X) \cong \Vect(B),
  \qquad
  \Vect_\Grd(X,A) \cong \Vect(B,B/J),
  \]
  which induce group isomorphisms
  \[
  \vbK^0_\Grd(X) \cong \K_{00}(B),
  \qquad
  \vbK^0_\Grd(X,A) \cong \K_{00}(B,B/J).
  \]
\end{theorem}

\begin{proof}
  This follows from
  Theorem~\ref{the:category_equivalence_equivariant_crossed}
  and Definitions \ref{def:vbK}, \ref{def:vbK_relative},
  \ref{def:Vect_B}, and~\ref{def:Vect_B_rel}.
\end{proof}

\subsection{Comparison with representable
  \texorpdfstring{$\K$}{K}-theory}
\label{sec:compare_vbK_RK}

With the same notation as in Theorem~\ref{the:vbK_via_crossed},
Corollaries \ref{cor:RK_via_crossed}
and~\ref{cor:RK_support_via_crossed} yield
\[
\RK^0_\Grd(X) \cong \K_0(B),
\qquad
\RK^0_\Grd(X,A) \cong \K_0(J).
\]
Thus it remains to compare the ``na\"\i ve'' \(\K\)\nb-theory
groups \(\K_{00}(B)\) and \(\K_{00}(B,B/J)\) with the correct
ones.  It is easy to construct natural transformations
\begin{align*}
  \vbK^0_\Grd(X) \cong \K_{00}(B) &\to \K_0(B) \cong \RK^0_\Grd(X),\\
  \vbK^0_\Grd(X,A) \cong \K_{00}(B,B/J) &\to \K_0(J) \cong
  \RK^0_\Grd(X,A).
\end{align*}

We briefly describe how the map
\(\vbK^0_\Grd(X,A)\to\RK^0_\Grd(X,A)\) acts on cycles for the
original definitions.  Let \((V^+,V^-,\varphi)\) be a cycle for
\(\vbK^0_\Grd(X,A)\), that is, \(V^+\) and~\(V^-\) are
\(\Grd\)\nb-equivariant vector bundles on~\(X\) and~\(\varphi\)
is an isomorphism \(V^+|_A\xrightarrow{\cong} V^-|_A\).  We
equip~\(V^\pm\) with \(\Grd\)\nb-invariant inner products, so
that their spaces of \(\CONT_0\)\nb-sections form
\(\Grd\)\nb-equivariant Hilbert modules over~\(\CONT_0(X)\).

We extend~\(\varphi\) to a \(\Grd\)\nb-equivariant, contractive
vector bundle map \(\cl{\varphi}\colon V^+\to V^-\): first use
the Tietze Extension Theorem to get a non-equivariant
extension; then use the properness of the action to replace
this by a \(\Grd\)\nb-equivariant function, as in the
definition of~\(T^\Grd\) before Proposition 6.24
in~\cite{Tu:Novikov}.  Since~\(\cl{\varphi}\) is invertible
on~\(A\), it remains invertible in a neighbourhood of~\(A\).
Using a partition of unity, we can adjust the Hermitian
structure so that~\(\cl{\varphi}\) is unitary in a
neighbourhood of~\(A\).

Hence we may also define \(\vbK^0_\Grd(X,A)\) using triples
\((V^+,V^-,\varphi)\), where \(V^+\) and~\(V^-\) are
\(\Grd\)\nb-equivariant Hermitian vector bundles on~\(X\) and
\(\varphi\colon V^+\to V^-\) is a \(\Grd\)\nb-equivariant
operator with \(\norm{\varphi}\le 1\) that is unitary in a
\(\Grd\)\nb-invariant neighbourhood of~\(A\).

Let \((V^+,V^-,\varphi)\) be a cycle for \(\vbK^0_\Grd(X,A)\)
of this new kind.  Then the spaces \(\Gamma_0(X\setminus
A,V^+)\) and \(\Gamma_0(X\setminus A,V^-)\) of
\(\CONT_0\)\nb-sections are \(\Grd\)\nb-equivariant Hilbert
modules over \(\CONT_0(X\setminus A)\), and~\(\varphi\) defines
an adjointable operator~\(M_\varphi\) between them.  The triple
\[
(\Gamma_0(X\setminus A,V^+),\Gamma_0(X\setminus A,V^-),
M_\varphi)
\]
is a cycle for \(\RK^0_{\Grd,X}(X\setminus A)\defeq
\KK^{\Grd\ltimes X}\bigl(\CONT_0(X),\CONT_0(X\setminus
A)\bigr)\) because the functions \(\ID_{V^+}-\varphi^*\varphi\)
and \(\ID_{V^-}-\varphi\varphi^*\) vanish outside an
\(X\)\nb-compact subset of \(X\setminus A\).  Notice that a
subset of \(X\setminus A\) is \(X\)\nb-compact if and only if
it is closed.

To get a map \(X\setminus A\to\Fred_\Grd\) from
\((V^+,V^-,\varphi)\), we add the degenerate cycle
corresponding to the identity operator on \(L^2(\Grd\ltimes
X\setminus A)^\infty\) and use the Equivariant Stabilisation
Theorem
\[
\Gamma_0(X\setminus A,V^+) \oplus L^2(\Grd\ltimes X\setminus A)^\infty
\cong L^2(\Grd\ltimes X\setminus A)^\infty.
\]
The resulting map \(U\colon X\setminus A\to\Fred_\Grd\) is
unitary where~\(\varphi\) is unitary, so that its support is
still \(X\)\nb-compact.

\smallskip

The natural transformation
\(\vbK^0_\Grd(X,A)\to\RK^0_\Grd(X,A)\) also yields one from
\(\vbK^{-n}_\Grd(X,A)\) to \(\RK^{-n}_\Grd(X,A)\) for all
\(n\in\N\), using~\eqref{eq:vbK_odd} and the natural
isomorphisms
\begin{equation}
  \label{eq:RK_odd}
  \RK^{-n}_\Grd(X,A) \cong
  \RK^0_\Grd(X\times\R^n,A\times\R^n) \cong
  \RK^0_\Grd(X\times\Sphere^n,X\times\{1\}\cup A\times\Sphere^n),
\end{equation}
where we use the excision property~\eqref{eq:RKG_excision}.

Now we have constructed a natural transformation
\(\vbK^*_\Grd(X,A)\to \RK^*_\Grd(X,A)\).  This is not an
isomorphism in general, as shown by several counterexamples
\cites{Phillips:Equivariant_Kbook, Lueck-Oliver:Completion,
  Sauer:K-theory}.  The following counterexample by Juliane
Sauer is particularly simple.

\begin{example}[see~\cite{Sauer:K-theory}]
  \label{exa:RK_not_vb}
  Consider the compact group \(K\defeq \prod_{n\in\Z} \Z/2\)
  and let~\(\alpha\) be the automorphism on~\(K\) that shifts
  the copies of~\(\Z/2\) to the left.  Form the semi-direct
  product group \(G\defeq \Z\ltimes_\alpha K\); this is a
  totally disconnected group.  Let~\(K\) act trivially
  and~\(\Z\) by translation on~\(\R\).  This combines to a
  proper, cocompact action of~\(G\) on~\(\R\).

  It is shown in~\cite{Sauer:K-theory} that the map
  \[
  \vbK^1_G(\R,\Z) \to \RK^1_G(\R,\Z)
  \]
  is not surjective.  The idea is the following.  On the right
  hand side, we compute
  \[
  \RK^1_G(\R,\Z) \cong \K^1_G(\R,\Z)
  \cong \K^1\bigl(G\ltimes \CONT_0(\R\setminus\Z)\bigr)
  \cong \K^1\bigl(\CONT_0\bigl((0,1)\bigr)\otimes C^*(K)\bigr)
  \cong \operatorname{Rep}(K).
  \]
  But any \(G\)\nb-equivariant vector bundle on
  \(\R\times\Sphere^1\) must carry the trivial representation
  of~\(K\) because this is the only representation that is
  fixed by the automorphism~\(\alpha\).  Hence
  \(\vbK^1_G(\R,\Z) = \vbK^1_\Z(\R,\Z)\cong\Z\) is much smaller
  than \(\RK^1_G(\R,\Z)\).

  A similar computation shows that the map
  \(\vbK^0_G(\R\times\Sphere^1)\to \RK^0_G(\R\times\Sphere^1)\)
  is not surjective, providing a counterexample for the
  absolute theories.

  By the way, since~\(\Z\) acts freely and properly on~\(\R\),
  the crossed product groupoid \(G\ltimes\R\) is Morita
  equivalent to a locally trivial bundle of compact groups on
  \(\Sphere^1 = \R/\Z\) whose fibre is~\(K\) everywhere.  More
  precisely, we consider the trivial bundle of groups
  \([0,1]\times K\) and identify \((0,g)\sim
  \bigl(1,\alpha(g)\bigr)\) for all \(g\in K\) to get a locally
  trivial bundle on~\(\Sphere^1\).  Since all our theories are
  evidently Morita invariant, we also get counterexamples for
  this groupoid, which is a locally trivial bundle of compact
  groups with compact base.
\end{example}

Several authors have established \(\vbK^*_G(X,A)\cong
\RK^*_G(X,A)\) for some classes of cocompact group actions.  Of
course, if~\(X\) is \(G\)\nb-compact then we can replace
\(\RK^*_G(X,A)\) by \(\K^*_G(X,A)\cong \K^*_G(X\setminus A)\).
If~\(G\) is almost connected or a matrix group, then Chris
Phillips (\cite{Phillips:Equivariant_K2}) proves an isomorphism
\(\vbK^*_G(X,A)\cong \K^*_G(X\setminus A)\) for all
\(G\)\nb-compact proper \(G\)\nb-spaces~\(X\).  Wolfgang L\"uck
and Bob Oliver (\cite{Lueck-Oliver:Completion}) show this
if~\(G\) is discrete and \((X,A)\) is a finite
\(G\)\nb-CW-pair.  The latter result is extended by Juliane
Sauer (\cite{Sauer:K-theory}) to totally disconnected groups
that are projective limits of discrete groups.  These results
are contained in Theorem~\ref{the:apprid_for_group} below.

Even without a group action, there are non-compact spaces for
which vector bundles do not generate the representable
\(\K\)\nb-theory.  It is known to experts that \(\vbK^*(X)\cong
\RK^*(X)\) if~\(X\) is, say, a finite-dimensional CW-complex;
more generally, this holds for paracompact Hausdorff spaces
with finite covering dimension.  But we could not find a proof
of this in the literature.  We do not consider equivariant
versions of this result here.

\section{Approximate units of projections}
\label{sec:cocompact_vbK}

\begin{definition}
  \label{def:appridpro}
  A (separable) \sigCstar{}algebra~\(B\) has an
  \emph{approximate unit of projections} if there is a sequence
  of projections \((p_n)\) in~\(B\) with \(\lim_{n\to\infty}
  p_n\cdot b\cdot p_n= b\) for all \(b\in B\).
\end{definition}

Since a sequence in \(B=\varprojlim B_n\) converges if and only
if its image in~\(B_n\) converges for all \(n\in\N\), a
sequence of projections \((p_n)\) in \(B=\varprojlim B_n\) is
an approximate unit in~\(B\) if and only if its image
in~\(B_n\) is an approximate unit for each \(n\in\N\).  But it
does not suffice merely to assume that all~\(B_n\) have
approximate units of projections because it is not clear
whether projections in~\(B_n\) lift to projections in~\(B\).

If~\(B\) has an approximate unit of projections, so
has~\(B_\Comp\): simply consider the sequence \((p_n\otimes
q_n)\) where \((p_n)\) and \((q_n)\) are approximate units of
projections in \(B\) and~\(\Comp\), respectively.  The
following example shows that the converse need not hold:

\begin{example}
  \label{exa:apprid_projection_stable}
  We describe a \Cstar{}algebra~\(B\) without an approximate
  unit of projections for which~\(B_\Comp\) has one.  Let
  \[
  B \defeq \bigl\{
  \left(\begin{smallmatrix}a&b\\c&d\end{smallmatrix}\right)
  \in \Mat_2\otimes \CONT\bigl([0,1]\bigr) \bigm|
  b(0)=c(0)=d(0)=0\bigr\}.
  \]
  This is the \Cstar{}algebra of a proper groupoid, namely, the
  groupoid associated to the equivalence relation on
  \[
  [0,1]\sqcup (0,1] = [0,1]\times\{0\} \cup (0,1]\times\{1\}
  \]
  that identifies \((t,0)\sim (t,1)\) for all \(t\in (0,1]\).
  Of course, this groupoid is Morita equivalent to \([0,1]\)
  with only identity morphisms.

  Correspondingly, \(B\) is Morita equivalent to \(C([0,1])\).
  Thus \(B_\Comp \cong C([0,1])_\Comp\).  This \Cstar{}algebra
  has an approximate unit of projections because \(C([0,1])\)
  is unital.  But if \(p\in B\) is a projection, then \(p(0)\)
  must have rank~\(1\), so that~\(p\) has rank~\(1\) on all of
  \([0,1]\).  Hence~\(B\) itself contains no approximate unit
  of projections.
\end{example}

Since our notions should be Morita invariant, it is better to
require only an approximate unit of projections in~\(B_\Comp\),
not in~\(B\).

A local Banach algebra~\(A\) for which \(\Mat_\infty(A)\) has
an approximate unit of projections is called stably unital in
\cite{Blackadar:K-theory}*{Definition 5.5.4}.  We avoid this
name for two reasons: first, \sigCstar{}algebras are local
Banach algebras if and only if they are \Cstar{}algebras;
secondly, we also consider \sigCstar{}algebras that are stably
isomorphic to unital \sigCstar{}algebras, and these also
deserve to be called stably unital.

\begin{proposition}
  \label{pro:apprid_K00}
  Let~\(B\) be a \Cstar{}algebra for which~\(B_\Comp\) has an
  approximate unit of projections.  Then \(\K_{00}(B) \cong
  \K_0(B)\).  If \(J\subseteq B\) is a closed ideal, then
  \(\K_{00}(B,B/J)\cong \K_0(J)\).
\end{proposition}

\begin{proof}
  The first assertion is \cite{Blackadar:K-theory}*{Proposition
    5.5.5}.  The second result is well-known if~\(B\) is unital
  (see \cite{Blackadar:K-theory}*{Theorem 5.4.2}).  We explain
  how to reduce the general case to this special case.

  Stabilising \(J\) and~\(B\), we may assume that~\(B\) itself
  has an approximate unit of projections \((p_n)_{n\in\N}\).
  The extension \(J\into B\prto B/J\) is the inductive limit in
  the category of \Cstar{}algebras of the extensions
  \[
  p_nJp_n \into p_nBp_n \prto \pi(p_n)(B/J)\pi(p_n).
  \]
  Since~\(\K_0\) is continuous for such inductive limits, and
  since the same holds for our relative theory
  \(\K_{00}(B,B/J)\), it suffices to prove
  \[
  \K_0(p_nJp_n) =
  \K_{00}\bigl(p_nBp_n,\pi(p_n)(B/J)\pi(p_n)\bigr)
  \]
  for all \(n\in\N\).  This reduces the assertion to the known
  case where~\(B\) is unital.
\end{proof}

Proposition~\ref{pro:apprid_K00} fails for \sigCstar{}algebras,
in general.  The main source of problems is that the set of
invertible elements in a \sigCstar{}algebra is not open.

\begin{theorem}
  \label{the:apprid_vbK}
  Let~\(\Grd\) be a second countable, locally compact,
  Hausdorff groupoid with Haar system and let~\(X\) be a
  proper, \(\Grd\)\nb-compact, second countable
  \(\Grd\)\nb-space.  Suppose that \(C^*(\Grd\ltimes
  X)_\Comp\) has an approximate unit of projections.

  Then the canonical map
  \[
  \vbK^*_\Grd(X,A) \to \RK^*_\Grd(X,A) \cong \K^*_\Grd(X,A)
  \cong \K^*_\Grd(X\setminus A)
  \]
  is an isomorphism for all closed \(\Grd\)\nb-invariant
  subsets \(A\subseteq X\).

  Let~\(Y\) be a \(\Grd\ltimes X\)\nb-space and let~\(S_X\)
  be the directed set of closed, \(\Grd\)\nb-invariant,
  \(X\)\nb-compact subsets of~\(Y\).  For \(A\in S_X\), let
  \(\partial A\) be the boundary of~\(A\) as a subset of~\(Y\).
  Then
  \[
  \RK^*_{\Grd,X}(Y) \cong
  \varinjlim_{A\in S_X} \vbK^*_\Grd(A,\partial A).
  \]
\end{theorem}

\begin{proof}
  Define \(B\) and~\(J\) as in
  Theorem~\ref{the:vbK_via_crossed}.  Thus
  \[
  \vbK^0_\Grd(X) \cong \K_{00}(B),
  \qquad
  \vbK^0_\Grd(X,A) \cong \K_{00}(B,B/J).
  \]
  Corollary~\ref{cor:RK_via_crossed} yields \(\RK^0_\Grd(X)
  \cong \K_0(B)\), and
  Corollary~\ref{cor:RK_support_via_crossed} and the definition
  of the relative theory in~\eqref{eq:RK_G_definition} yield
  \(\RK^0_\Grd(X,A) \cong \K_0(J)\).  Hence the first two
  assertions of the theorem follow from
  Proposition~\ref{pro:apprid_K00}.

  To prove the last assertion, we claim that there are natural
  isomorphisms
  \[
  \RK^*_{\Grd,X}(Y)
  \cong \varinjlim_{A\in S_X} \RK^*_\Grd(Y,Y\setminus A^\circ)
  \cong \varinjlim_{A\in S_X} \RK^*_\Grd(A,\partial A)
  \cong \varinjlim_{A\in S_X} \vbK^*_\Grd(A,\partial A).
  \]
  Here~\(A^\circ\) denotes the interior of~\(A\), so that
  \(Y\setminus A^\circ\) is closed in~\(Y\) for all \(A\in
  S_X\).  The first isomorphism follows from
  Theorem~\ref{the:RK_support_continuous}.  The second
  isomorphism is the excision
  isomorphism~\eqref{eq:RKG_excision} for \(\RK^*_\Grd\).  For
  the third one, we use that \(\Grd\ltimes \CONT_0(A)_\Comp\)
  has an approximate unit of projections as well.  This is true
  because the proper map \(A\to X\) induces an essential
  \Star{}homomorphism \(\Grd\ltimes \CONT_0(X)_\Comp \to
  \Grd\ltimes \CONT_0(A)_\Comp\), which maps the approximate
  unit of projections in \(\Grd\ltimes \CONT_0(X)_\Comp\) to
  one in \(\Grd\ltimes \CONT_0(A)_\Comp\).
\end{proof}

The inductive limit \(\varinjlim_{A\in S_X}
\vbK^*_\Grd(A,\partial A)\) implicitly uses maps
\[
\vbK^*_\Grd(A,\partial A) \to \vbK^*_\Grd(B,\partial B)
\]
for \(A,B\in S_X\) with \(A\subseteq B\).  To construct these,
we need non-trivial excision isomorphisms
\[
\vbK^*_\Grd(A,\partial A) \cong \vbK^*_\Grd(B,B\setminus A^\circ).
\]
We have excision here because both sides are naturally isomorphic to
\(\RK^*_\Grd(A,\partial A)\) and \(\RK^*_\Grd(B,B\setminus A^\circ)\),
respectively, which do satisfy excision by~\eqref{eq:RKG_excision}.

It will be useful for later to formalise an idea in the proof of
Theorem~\ref{the:apprid_vbK}:

\begin{proposition}
  \label{pro:inherit_apprid}
  Let \(f\colon X\to Y\) be a \(\Grd\)\nb-equivariant
  continuous map between two proper \(\Grd\)\nb-spaces.  If
  \(\sC(\Grd\ltimes Y)_\Comp\) has an approximate unit of
  projections, then so has \(\sC(\Grd\ltimes X)_\Comp\).

  In particular, let \(\EG\Grd\) be a universal proper
  \(\Grd\)\nb-space.  If \(\sC(\Grd\ltimes \EG\Grd)_\Comp\)
  has an approximate unit of projections, then
  \(\sC(\Grd\ltimes X)_\Comp\) has an approximate unit of
  projections for all proper \(\Grd\)\nb-spaces~\(X\).
\end{proposition}

\begin{proof}
  Choose exhausting sequences \((X_n)\) and \((Y_n)\) of
  \(\Grd\)\nb-compact subsets as in
  Theorem~\ref{the:RKK_via_crossed}.  The map \(f\colon X\to
  Y\) restricts to proper maps \(X_n\to Y_n\), no matter
  whether~\(f\) itself is proper.  This is because \(X_n\)
  and~\(Y_n\) are \(\Grd\)\nb-compact and proper, so that any
  continuous \(\Grd\)\nb-equivariant map between them is
  proper.  Hence we get induced \Star{}homomorphisms
  \(\Grd\ltimes\CONT_0(Y_n)\to\Grd\ltimes\CONT_0(X_n)\), which
  are essential.  The resulting \Star{}homomorphism
  \[
  \sC(\Grd\ltimes Y)_\Comp \to \sC(\Grd\ltimes X)_\Comp
  \]
  maps an approximate unit of projections for~\(Y\) to one
  for~\(X\), yielding the first assertion.  The last assertion
  follows because any proper \(\Grd\)\nb-space maps to
  \(\EG\Grd\).
\end{proof}

\subsection{General existence criteria}
\label{sec:apprid_exist_general}

Theorem~\ref{the:apprid_vbK} motivates us to search for
approximate units of projections in stable \Cstar{}algebras and
\sigCstar{}algebras.  We will see that there is an approximate
unit of projections in the stabilisation if there are enough
projections in a weaker sense.  Various counterexamples suggest
that we cannot do much better here.

To simplify our notation, we consider only stable
\sigCstar{}algebras from now on, that is, we assume \(B\cong
B_\Comp\).  This can be achieved by replacing~\(B\)
by~\(B_\Comp\).

\begin{definition}
  \label{def:prid}
  Let \(\Prid(B)\) denote the primitive ideal space of a
  \Cstar{}algebra~\(B\), equipped with the hull-kernel
  topology.  For a \sigCstar{}algebra \(B=\varprojlim B_n\),
  the projections \(B_{n+1}\prto B_n\) yield an inductive
  system of closed embeddings of topological spaces
  \(\Prid(B_n)\to \Prid(B_{n+1})\).  We let \(\Prid(B)\) be the
  direct limit of this system or, less formally, \(\Prid(B)
  \defeq \bigcup_{n=1}^\infty \Prid(B_n)\).
\end{definition}

For each primitive ideal \(\mathfrak{p}\in\Prid(B_n)\subseteq
\Prid(B)\), we get a \Cstar{}algebra quotient \(B/\mathfrak{p}
\defeq B_n/\mathfrak{p}\), which does not depend on the choice
of~\(n\).

\begin{definition}
  \label{def:sigma_unital}
  A \sigCstar{}algebra \(B=\varprojlim B_n\) is called
  \emph{\(\sigma\)\nb-unital} if it contains a strictly
  positive element \(h\in B\).
\end{definition}

It is easy to see that \(B=\varprojlim B_n\) is
\(\sigma\)\nb-unital if and only if~\(B_n\) is
\(\sigma\)\nb-unital for all \(n\in\N\).  This holds, for
instance, if all~\(B_n\) are separable.

\begin{theorem}
  \label{the:apprid_exists}
  A stable, \(\sigma\)\nb-unital \sigCstar{}algebra~\(B\)
  contains an approximate unit of projections if and only if
  for each primitive ideal~\(\mathfrak{p}\) in~\(B\) there is a
  projection~\(q\) in~\(B\) that does \emph{not} belong
  to~\(\mathfrak{p}\).
\end{theorem}

\begin{proof}
  Suppose first that \(B\cong B_\Comp\) contains an approximate
  unit of projections \((e_n)_{n\in\N}\).  If \(\mathfrak{p}\in
  B\), then the images of~\(e_n\) in \(B/\mathfrak{p}\) form an
  approximate unit of projections as well.  Hence \(e_n\notin
  \mathfrak{p}\) for some \(n\in\N\).

  The converse direction is more interesting.  For a projection
  \(q\in B\), let~\(I_q\) be the closed \Star{}ideal generated
  by~\(q\), that is, the closed linear span of \(aqb\) for
  \(a,b\in B\).  The ideal~\(I_q\) is Morita equivalent to the
  corner \(qBq\) via the imprimitivity bimodule \(qB\)
  (see~\cite{Brown:Stable_isomorphism}).  Thus each
  ideal~\(I_q\) is Morita equivalent to a unital
  \Cstar{}algebra.  We identify \(\Prid(I_q)\) with an open
  subset in \(\Prid(B)\) in the usual way.  We have
  \(\mathfrak{p}\in \Prid(I_q)\) if and only if
  \(q\notin\mathfrak{p}\).  Hence the subsets \(\Prid(I_q)\)
  for projections \(q\in B\) form an open covering of
  \(\Prid(B)\).  More precisely, they form open coverings of
  \(\Prid(B_n)\) for each \(n\in\N\).

  Let \(h\in B\) be strictly positive.  The subsets
  \[
  X_n \defeq
  \{\mathfrak{p}\in\Prid(B) \mid
  \norm{h}_{B/\mathfrak{p}}\ge \nicefrac{1}{n}\}
  \]
  for \(n\in\N\) have quasi-compact intersection with
  \(\Prid(B_k)\) for each \(k\in\N\).  They form an increasing
  sequence with \(\bigcup_{n\in\N} X_n=\Prid(B)\) because
  \(h\in B\) is strictly positive.  Let \(n\in\N\).  Since
  \(X_n\cap \Prid(B_n)\) is quasi-compact, we can find
  projections \(q_1,\dotsc,q_j\in B\) with
  \[
  X_n\cap \Prid(B_n)\subseteq
  \Prid(I_{q_1})\cup\dotsb \cup \Prid(I_{q_j}).
  \]
  The projection \(q_1\oplus\dotsb\oplus q_j\in \Mat_j(B)\) is
  Murray--von Neumann equivalent to a projection~\(e_n\)
  in~\(B\) because~\(B\) is stable.  It is easy to see that
  Murray--von Neumann equivalent projections generate the same
  ideal.  Hence the primitive ideal space of \(I_n\defeq
  I_{e_n}\) contains \(X_n\cap\Prid(B_n)\).  Since \(\bigcup
  X_n=\Prid(B)\), we get \(\bigcup \Prid(I_n)=\Prid(B)\).  This
  implies that \(\bigcup I_n\) is dense in~\(B\): a subset is
  dense in~\(B\) if and only if its image in~\(B_k\) is dense
  for all \(k\in\N\), and open subsets of \(\Prid(B_k)\) are in
  bijection with ideals in~\(B_k\).

  Each of the ideals~\(I_n\) is Morita equivalent to a unital
  \sigCstar{}algebra, namely, \(e_nBe_n\).  Therefore,
  \((I_n)_\Comp\) contains an approximate unit of projections
  \((e_{n,k})_{k\in\N}\).  Since \(\bigcup (I_n)_\Comp\) is
  dense in~\(B_\Comp\), \(e'_n\defeq e_{n,k_n}\) for a suitable
  function \(n\mapsto k_n\) is an approximate unit of
  projections in~\(B_\Comp\).
\end{proof}

\begin{definition}
  \label{def:full_projection}
  A set of projections~\(S\) is called \emph{full} if for any
  primitive ideal~\(\mathfrak{p}\) in~\(B\), there is \(q\in
  S\) that does \emph{not} belong to~\(\mathfrak{p}\).  A
  single projection~\(q\) is called full if \(\{q\}\) is.
\end{definition}

Thus Theorem~\ref{the:apprid_exists} asserts that a
\(\sigma\)\nb-unital stable \sigCstar{}algebra has an
approximate unit of projections if and only if it contains a
full set of projections.

\begin{corollary}
  \label{cor:apprid_compact}
  The following are equivalent for a \(\sigma\)\nb-unital stable
  \Cstar{}algebra~\(B\):
  \begin{enumerate}[label=\textup{(\arabic{*})}]
  \item \(B\) contains a full projection;
  \item \(B\) is Morita equivalent to a unital algebra;
  \item \(B\) contains an approximate unit of projections and
    \(\Prid(B)\) is quasi-compact;
  \item \(B\) contains a full set of projections and
    \(\Prid(B)\) is quasi-compact.
  \end{enumerate}
  Statements \textup{(1)} and~\textup{(2)} remain equivalent
  for a \(\sigma\)\nb-unital stable \sigCstar{}algebra~\(B\).
\end{corollary}

\begin{proof}
  (1)\(\iff\)(2) for \sigCstar{}algebras: If \(p\in B\) is a
  full projection, then the ideal~\(I_p\) generated by~\(p\) is
  all of~\(B\), and Morita equivalent to the unital algebra
  \(pBp\).  Conversely, if~\(B\) is Morita equivalent to a
  unital \sigCstar{}algebra~\(A\), then \(B\cong B_\Comp\cong
  A_\Comp\), and the class of the unit in~\(A\) is a full
  projection in~\(A_\Comp\).

  Now we assume that~\(B\) is a \Cstar{}algebra.

  (2)\(\Longrightarrow\)(3): since~\(B\) is stable, (2) implies
  \(B\cong A_\Comp\) for a unital \Cstar{}algebra~\(A\).  Then
  \(1_A\otimes e_n\) is an approximate unit of projections
  in~\(A\), where \((e_n)\) is one in~\(\Comp\), and \(\Prid(A)
  \cong \Prid(B)\) is quasi-compact.

  (3)\(\Longrightarrow\)(4) is trivial.

  (4)\(\Longrightarrow\)(1): This follows from the proof of
  Theorem~\ref{the:apprid_exists}.  Since \(\Prid(B)\) is
  quasi-compact, we have \(\Prid(B)=\Prid(I_n)\) for some
  \(n\in\N\).  Hence the projection~\(q_n\) that
  generates~\(I_n\) is full.
\end{proof}

\begin{corollary}
  \label{cor:apprid_prid}
  Let~\(B\) be a \(\sigma\)\nb-unital \Cstar{}algebra.
  If~\(B_\Comp\) contains an approximate unit of projections,
  then there is an increasing sequence of quasi-compact open
  subsets \(U_n\subseteq \Prid(B)\) with \(\bigcup
  U_n=\Prid(B)\).
\end{corollary}

\begin{proof}
  This follows from the proof of
  Theorem~\ref{the:apprid_exists}.  The subsets \(\Prid(I_n)\)
  are open and quasi-compact because each~\(I_n\) is Morita
  equivalent to a unital \Cstar{}algebra.
\end{proof}

\begin{example}
  \label{exa:circle_group_bundle}
  Consider again the situation of Example~\ref{exa:RK_not_vb},
  that is, let~\(\Grd\) be (Morita equivalent) to the locally
  trivial bundle of compact groups on~\(\Sphere^1\) with
  fibre~\(K\), twisted by the automorphism~\(\alpha\).  The
  primitive ideal space~\(P\) of \(C^*(\Grd)\) is Hausdorff in
  this case, and the map \(P\to \Sphere^1\) is a covering map
  with fibre \(\Prid(C^*K) = \widehat{K}\).  The holonomy of
  this covering is given by the action of the
  automorphism~\(\alpha\) on~\(\widehat{K}\).  More explicitly,
  \[
  \widehat{K} = \bigoplus_{n\in\Z} \widehat{\Z/2}
  \cong \bigoplus_{n\in\Z} \Z/2,
  \]
  with~\(\alpha\) acting by translation.  This action is free
  on \(\widehat{K}\setminus\{\tau\}\), where~\(\tau\) denotes
  the trivial representation.  Thus \(\Prid(C^*\Grd)\) contains
  several copies of~\(\R\) --~one for each orbit of~\(\alpha\)
  on \(\widehat{K}\setminus\{\tau\}\)~-- and one copy
  of~\(\Sphere^1\).  Since~\(\R\) contains no quasi-compact
  open subsets, the necessary condition of
  Corollary~\ref{cor:apprid_prid} is violated in this case.
\end{example}

\begin{example}
  \label{exa:continuous_trace}
  The necessary condition in Corollary~\ref{cor:apprid_prid} is
  not sufficient, even if we restrict attention to, say, type~I
  \Cstar{}algebras.  Counterexamples come from continuous trace
  \Cstar{}algebras (see
  also~\cite{Tu-Xu-Laurent-Gengoux:Twisted_K}).  Let~\(B\) be a
  continuous trace \Cstar{}algebra with \emph{connected}
  spectrum~\(X\).  Stable isomorphism classes of continuous
  trace \Cstar{}algebras with spectrum~\(X\) are classified by
  their Dixmier--Douady invariant in the \v{C}ech cohomology
  \(H^3(X,\Z)\) (see \cite{Cuntz-Meyer-Rosenberg}*{Theorem
    9.9}).  If~\(p\) is a projection in a stable continuous
  trace \Cstar{}algebra~\(B\) over~\(X\), then the rank
  of~\(p\) is a constant function on~\(X\).  Thus \(pBp\) is a
  locally trivial bundle of finite-dimensional matrix algebras.
  But this implies that the Dixmier--Douady invariant in
  \(H^3(X,\Z)\) is torsion (see
  \cite{Cuntz-Meyer-Rosenberg}*{Theorem 9.13}).

  Hence a continuous trace \Cstar{}algebra with non-torsion
  Dixmier--Douady invariant contains no projections, even
  stably.  In contrast, a continuous trace \Cstar{}algebra with
  torsion Dixmier--Douady invariant is Morita equivalent to a
  unital \Cstar{}algebra and hence contains an approximate unit
  of projections in its stabilisation by
  Corollary~\ref{cor:apprid_compact}.
\end{example}

Example~\ref{exa:continuous_trace} shows that the existence of
an approximate unit of projections in the stabilisation is a
very subtle global question.  This explains why our criteria
for groupoid \Cstar{}algebras require the existence of some
equivariant vector bundles to begin with.

\subsection{Application to groupoids}
\label{sec:apprid_groupoid}

\begin{theorem}
  \label{the:groupoid_apprid}
  Let~\(\Grd\) be a groupoid and let~\(X\) be a proper
  \(\Grd\)\nb-space.  The \sigCstar{}algebra \(\sC(\Grd\ltimes
  X)_\Comp\) has an approximate unit of projections if and only
  if, for each \(x\in X\) and each irreducible representation
  \(\varrho\colon \Grd^x_x\to\Mat_n(\C)\) of its
  stabiliser~\(\Grd_x^x\), there is a \(\Grd\)\nb-equivariant
  vector bundle on~\(X\) whose fibre at~\(x\) contains the
  representation~\(\varrho\).

  The \sigCstar{}algebra \(\sC(\Grd\ltimes X)_\Comp\) is Morita
  equivalent to a unital \sigCstar{}algebra if and only if
  there is a single \(\Grd\)\nb-equivariant vector bundle
  on~\(X\) whose fibre at~\(x\) contains all irreducible
  representations of~\(\Grd_x^x\) for all \(x\in X\).
\end{theorem}

\begin{proof}
  First we recall how to compute the primitive ideal space of
  \(\sC(\Grd\ltimes X)\).  It suffices to do this in the
  cocompact case where \(\sC(\Grd\ltimes X)\) is a
  \Cstar{}algebra because \(\Prid\bigl(\sC(\Grd\ltimes X)\bigr)
  = \bigcup \Prid\bigl(C^*(\Grd\ltimes X_n)\bigr)\);
  here~\(X_n\) is an increasing sequence of cocompact
  \(\Grd\)\nb-invariant subsets with \(\bigcup X_n=X\).

  Any irreducible \Star{}representation of \(C^*(\Grd\ltimes
  X_n)\) is carried by a single orbit because the central
  subalgebra \(\CONT_0(\Grd\backslash X_n)\) of
  \(\Mult\bigl(C^*(\Grd\ltimes X_n)\bigr)\) must act by a
  character.  The irreducible representations carried by the
  orbit \(\Grd\cdot x\) correspond bijectively to the
  irreducible representations of the restricted groupoid
  \(\Grd\ltimes (\Grd\cdot x)\), which is Morita equivalent to
  the stabiliser~\(\Grd_x^x\) of~\(x\) because it is proper and
  transitive.  Thus irreducible representations correspond
  bijectively to \(\bigsqcup_{x\in R} \widehat{\Grd_x^x}\),
  where~\(R\) is a set of representatives for the
  \(\Grd\)\nb-orbits in~\(X\).  It follows easily that each
  irreducible representation of \(C^*(\Grd\ltimes X_n)\) is
  completely continuous.  Hence the primitive ideal space agrees
  with the space of irreducible representations.

  Theorem~\ref{the:category_equivalence_equivariant_crossed}
  shows that equivalence classes of idempotents in
  \(\sC(\Grd\ltimes X)_\Comp\) correspond bijectively to
  isomorphism classes of \(\Grd\)\nb-equivariant vector bundles
  on~\(X\).  The projection associated to a
  \(\Grd\)\nb-equivariant vector bundle~\(V\) vanishes at a
  representation \(\varrho\in\widehat{\Grd_x^x}\) if and only
  if the fibre~\(V_x\), which is a representation
  of~\(\Grd_x^x\), does not contain~\(\varrho\).  Hence the
  assertion follows from Theorem~\ref{the:apprid_exists}.
  Furthermore, a single vector bundle contains all irreducible
  representations of stabilisers if and only if the
  corresponding projection in \(\sC(\Grd\ltimes X)_\Comp\) is
  full.  Such a projection exists if and only if
  \(\sC(\Grd\ltimes X)_\Comp\) is Morita equivalent to a unital
  \sigCstar{}algebra.
\end{proof}

Theorem~\ref{the:groupoid_apprid} seems disappointing at first
sight because we require the existence of some equivariant
vector bundles to begin with.  Examples
\ref{exa:circle_group_bundle} and~\ref{exa:continuous_trace}
suggest that we must construct some equivariant vector bundles
by hand.

The following theorem contains the three previously known
situations where equivariant \(\K\)\nb-theory for group actions
can be computed by equivariant vector bundles (see
\cites{Phillips:Equivariant_K2, Lueck-Oliver:Completion,
  Sauer:K-theory}).

\begin{theorem}
  \label{the:apprid_for_group}
  Let~\(G\) be a locally compact group and let~\(X\) be a
  proper \(G\)\nb-compact \(G\)\nb-space.  The
  \sigCstar{}algebra \(\sC(G\ltimes X)_\Comp\) has an
  approximate unit of projections in each of the following
  cases:
  \begin{enumerate}[label=\textup{(\arabic{*})}]
  \item if~\(G\) is a closed subgroup of an almost connected
    group~\(H\), that is, the component group \(\pi_0(H)\) is
    compact;

  \item if~\(G\) is discrete, \(G\backslash X\) has finite
    covering dimension, and all finite subgroups of~\(G\) have
    order at most~\(N\) for some \(N\in\N\);

  \item more generally, it suffices to assume that there is a
    decreasing sequence of compact, open, normal subgroups
    \((K_n)_{n\in\N}\) such that
    \begin{itemize}
    \item \(\bigcap K_n = \{1\}\);
    \item for each \(n\in\N\) there is \(N_n\in\N\) such that
      all finite subgroups of the discrete group \(G/K_n\) have
      at most~\(N_n\) elements;
    \item the orbit spaces \(K_n\backslash X\) have finite
      covering dimension for all \(n\in\N\).
    \end{itemize}
  \end{enumerate}
  Furthermore, in the second case, \(\sC(G\ltimes X)_\Comp\)
  contains a full projection, so that \(\sC(G\ltimes X)\) is
  Morita equivalent to a unital \sigCstar{}algebra.
\end{theorem}

\begin{proof}
  We treat case (1) first.  We induce the action of~\(G\) to an
  action \(H\times_G X\) of~\(H\), which is still cocompact.
  The groupoids \(G\ltimes X\) and \(H\ltimes (H\times_G X)\)
  are Morita equivalent, so that it makes no difference which
  one we study.  Let \(K\subseteq H\) be a maximal compact
  subgroup.  Then the quotient space \(H/K\) is a universal
  proper \(H\)\nb-space by~\cite{Abels:Slices}.  This model
  of \(\EG H\) is \(H\)\nb-compact, and \(C^*(H\ltimes H/K)\)
  is Morita equivalent to \(C^*(K)\).  Since the latter has an
  approximate unit of projections,
  Proposition~\ref{pro:inherit_apprid} yields one in
  \(\sC\bigl(H\ltimes (H\times_G X)\bigr)_\Comp\) and hence in
  \(\sC(G\ltimes X)_\Comp\).

  Now we turn to case (2).  For a discrete group~\(G\),
  Wolfgang L\"uck and Bob Oliver construct a full
  \(G\)\nb-equivariant vector bundle on any
  finite-dimensional proper \(G\)\nb-CW-complex whose
  isotropy groups have bounded order (see
  \cite{Lueck-Oliver:Completion}*{Corollary 2.7}).
  Proposition~\ref{pro:inherit_apprid} yields our assertion
  for~\(X\) if we can find a continuous \(G\)\nb-map
  from~\(X\) to such a space.  This is an easy application of
  partitions of unity.

  Any proper group action is locally induced, that is, there is
  an open \(G\)\nb-invariant covering~\(\mathcal{U}\) such
  that for each \(U\in\mathcal{U}\), there is a
  \(G\)\nb-equivariant homeomorphism \(U\cong G\times_H Y\)
  for some finite subgroup \(H\subseteq G\) and some
  \(H\)\nb-space~\(Y\) (see \cites{Abels:Universal,
    Chabert-Echterhoff-Meyer:Deux}).  Equivalently, there is a
  \(G\)\nb-equivariant map \(U\to G/H\).  Since the orbit
  space \(G\backslash X\) has finite covering dimension, we may
  assume that \(\mathcal{U}\) splits into \(n+1\) subsets
  \(\mathcal{U}_0,\dotsc,\mathcal{U}_n\) such that \(U\cap
  V=\emptyset\) for all \(U,V\in\mathcal{U}_j\) for the
  same~\(j\).  Let \((\varphi_U)_{U\in\mathcal{U}}\) be a
  partition of unity on~\(X\) by \(G\)\nb-invariant functions
  subordinate to this covering.

  Since each \(U\in\mathcal{U}\) comes with a map \(f\colon
  U\to G/H\) for some finite subgroup~\(H\), we can
  decompose~\(U\) into a disjoint union \(U=\bigsqcup_{gH\in
    G/H} f^{-1}(gH)\).  We let~\(\mathcal{U}'\) be the
  resulting open covering by these subsets.  It still has the
  same finite covering dimension, so that its nerve
  \(\abs{\mathcal{U}'}\) is a finite-dimensional simplicial
  complex.  The obvious action of~\(G\) on~\(\mathcal{U}'\)
  induces a proper simplicial action on \(\abs{\mathcal{U}'}\).
  Let \(\varphi_{f^{-1}(gH)} \defeq \varphi_U|_{f^{-1}(gH)}\).
  This is a \(G\)\nb-equivariant partition of unity on~\(X\).
  It induces a \(G\)\nb-equivariant continuous map \(X\to
  Y\).  The assumption on~\(G\) ensures that the isotropy
  groups in~\(\abs{\mathcal{U}'}\) have finite order.

  Thus \cite{Lueck-Oliver:Completion}*{Corollary 2.7} provides
  a full \(G\)\nb-equivariant vector bundle
  on~\(\abs{\mathcal{U'}}\), which pulls back to a full
  \(G\)\nb-equivariant vector bundle on~\(X\).  Now
  Theorem~\ref{the:groupoid_apprid} yields a full projection
  and an approximate unit of projections in \(\sC(G\ltimes
  X)_\Comp\).

  Finally, we reduce case (3) to case (2).  There are full
  \(G/K_n\)-equivariant vector bundles on \(K_n\backslash
  X\) for all \(n\in\N\) by (2).  Pull them back to
  \(G\)\nb-equivariant vector bundles~\(V_n\) on~\(X\).  For
  \(x\in X\), any representation of the stabiliser~\(G_x^x\)
  must be trivial on~\(K_n\) for sufficiently large~\(n\) and
  hence be contained in the fibre of~\(V_n\).  Hence the set of
  equivariant vector bundles~\(V_n\) is full.
\end{proof}

If Theorem~\ref{the:apprid_for_group} applies and the action
of~\(G\) on~\(X\) is cocompact, then
Theorem~\ref{the:apprid_vbK} shows that \(\RK^0_G(X)\) is
generated by vector bundles.

\begin{example}
  \label{exa:trivial_bundle_full}
  If a groupoid~\(\Grd\) acts freely on~\(X\), that is, all
  stabilisers are trivial, then the trivial
  \(1\)\nb-dimensional vector bundle with the obvious
  \(\Grd\)\nb-action is full, so that \(\sC(\Grd\ltimes X)\)
  is Morita equivalent to a unital \sigCstar{}algebra.
  If~\(\Grd\) already acts freely on its object space, that is,
  \(\Grd\) describes an equivalence relation, then any action
  of~\(\Grd\) is free, so that the above applies to all proper
  actions of~\(\Grd\).
\end{example}

\begin{example}
  \label{exa:orbifold_groupoid}
  Recall that orbifolds are described by effective proper
  \'etale Lie groupoids.
  \cite{Moerdijk-Pronk:Orbifolds}*{Theorem 4.1} asserts that
  any orbifold groupoid~\(\Grd\) is Morita equivalent to
  \(L\ltimes X\) for an action of a compact group~\(L\) on a
  smooth manifold~\(X\).  More precisely, we can let~\(X\) be
  the \(\Grd\)\nb-orbit space of the frame bundle
  of~\(\Grd^{(0)}\) and~\(L\) be the orthogonal group of
  appropriate dimension.

  Hence \(\sC(\Grd)_\Comp\) contains an approximate unit of
  projections, so that the equivariant \(\K\)\nb-theory of a
  cocompact orbifold groupoid is generated by equivariant
  vector bundles.
\end{example}

\section{Conclusion}
\label{sec:conclusion}

We have defined \(\K\)\nb-theory groups with several support
conditions that are equivariant with respect to proper actions
of groupoids.  While they are originally defined using
bivariant \(\K\)\nb-theory for \(C^*\)\nb-algebras, their
topological nature is made clear by alternative descriptions
using homotopy classes of equivariant maps to suitable spaces
of Fredholm operators.  Another description using
\(\sigma\)\nb-\(C^*\)-algebras is most convenient for
establishing the formal properties of these theories.

A difficult question is when equivariant \(\K\)\nb-theory can
be described by equivariant vector bundles or, more precisely,
whether the Grothendieck group of the monoid of equivariant
vector bundles on a space agrees with its equivariant
representable \(\K\)\nb-theory.  For cocompact actions, we have
found a useful criterion that allows to study this question in
examples, reducing it to the question: which stable
\(C^*\)\nb-algebras contain an approximate unit of projections?

New difficulties appear for actions that are not cocompact.
Here we should assume finite covering dimension of the orbit
space and the existence of a full equivariant vector bundle to
achieve anything.  So far, there appear to be no general
results in this case that involve non-trivial groups.  We plan
to remedy this in a forthcoming article.

% \affiliationone{Heath Emerson\\
%   Department of Mathematics and Statistics\\
%   University of Victoria\\
%   PO BOX 3045 STN CSC\\
%   Victoria, B.C.\\
%   Canada V8W 3P4
%   \email{hemerson@math.uvic.ca}}
% \affiliationtwo{Ralf Meyer\\
%   Mathematisches Institut\\
%   Georg-August-Universit\"at G\"ottingen\\
%   Bunsenstra{\ss}e 3--5\\
%   37073 G\"ottingen\\
%   Germany
%   \email{rameyer@uni-math.gwdg.de}}

\end{document}